\documentclass[11pt]{amsart}
\usepackage{a4wide}
\usepackage{comment}
\usepackage[dvipsnames]{xcolor}
\usepackage{mathrsfs}
\usepackage{amsmath}
\usepackage{amsthm}
\usepackage{amsfonts}
\usepackage{amssymb}
\usepackage{enumitem}
\usepackage[margin=1in,headsep=.2in,marginpar=0.5in]{geometry}
\usepackage{graphicx}
\usepackage{url}
\usepackage{tikz,pgfplots}
    \pgfplotsset{compat=1.15}
    \usetikzlibrary{arrows.meta}
    \usetikzlibrary{decorations.markings}
    \usetikzlibrary{calc}
    \usetikzlibrary{knots}
    \usetikzlibrary{external}
\usepackage{tikz-cd}
\usepackage{tabularray}
	\SetTblrInner{rowsep=0pt}
\usepackage{palatino}
\usepackage{overpic}
\usepackage{subcaption}
\usepackage[normalem]{ulem}
\usepackage[T1]{fontenc}
\usepackage{hyperref}
\usepackage[noabbrev,capitalize,nameinlink]{cleveref}
\newtheorem{theorem}{Theorem}[section]
\newtheorem{proposition}[theorem]{Proposition}
\newtheorem{lemma}[theorem]{Lemma}

\theoremstyle{definition}
\newtheorem{definition}[theorem]{Definition}

\newtheorem{question}[theorem]{Question}

\theoremstyle{remark}
\newtheorem{remark}[theorem]{Remark}
\numberwithin{equation}{section}

\newcommand{\Xne}{X\!\colon\!\!NE}
\newcommand{\Xnw}{X\!\colon\!\!NW}
\newcommand{\Xsw}{X\!\colon\!\!SW}
\newcommand{\Xse}{X\!\colon\!\!SE}
\newcommand{\One}{O\!\colon\!\!NE}
\newcommand{\Onw}{O\!\colon\!\!NW}
\newcommand{\Osw}{O\!\colon\!\!SW}
\newcommand{\Ose}{O\!\colon\!\!SE}

\DeclareMathOperator{\Fix}{Fix}
\DeclareMathOperator{\tb}{tb}
\DeclareMathOperator{\rot}{rot}

\definecolor{darkgreen}{RGB}{0,153,0}
\hypersetup{colorlinks,linkcolor={blue},citecolor={darkgreen},urlcolor={darkgreen}}

\newcommand{\jr}{\textcolor{blue}}

\begin{document}

\title{Size-preserving Cromwell equivalence for unknotted grid diagrams}

\author{Austin Christian}
\author{James Rea}
\author{Matthew Stark}
\address{California Polytechnic State University, San Luis Obispo, CA 93407}
\urladdr{https://sites.google.com/view/austin-christian}

\begin{abstract}
We demonstrate that every grid diagram representing the unknot is equivalent, via Cromwell moves which preserve the grid number, to a grid diagram with no crossings, sharpening a monotonic-simplification result of Dynnikov. Moreover, we enumerate the equivalence classes of grid diagrams for the unknot, considered up to size-preserving Cromwell moves. As an application, we generate all grid diagrams for the unknot through grid number $8$ and identify the first examples of weakly Cromwell hard grid diagrams for the unknot.
\end{abstract}

\maketitle

\tableofcontents

\section{Introduction}
	
	Grid diagrams have long been used to represent knots and links in $\mathbb{R}^3$ (cf. \cite{brunn1897uber,cromwell1995embedding,dynnikov2006arc}), and have been especially useful in the present century for computations of knot/link Floer homology (cf. \cite{manolescu2007combinatorial,manolescu2009combinatorial,manolescu2025grid}) and manipulating Legendrian knots and links (cf. \cite{chongchitmate2013atlas,petkova2026legendrian,lsuatlas,jubeir2024spectral}).
	
	\begin{figure}[ht]
	    \centering
	    \input{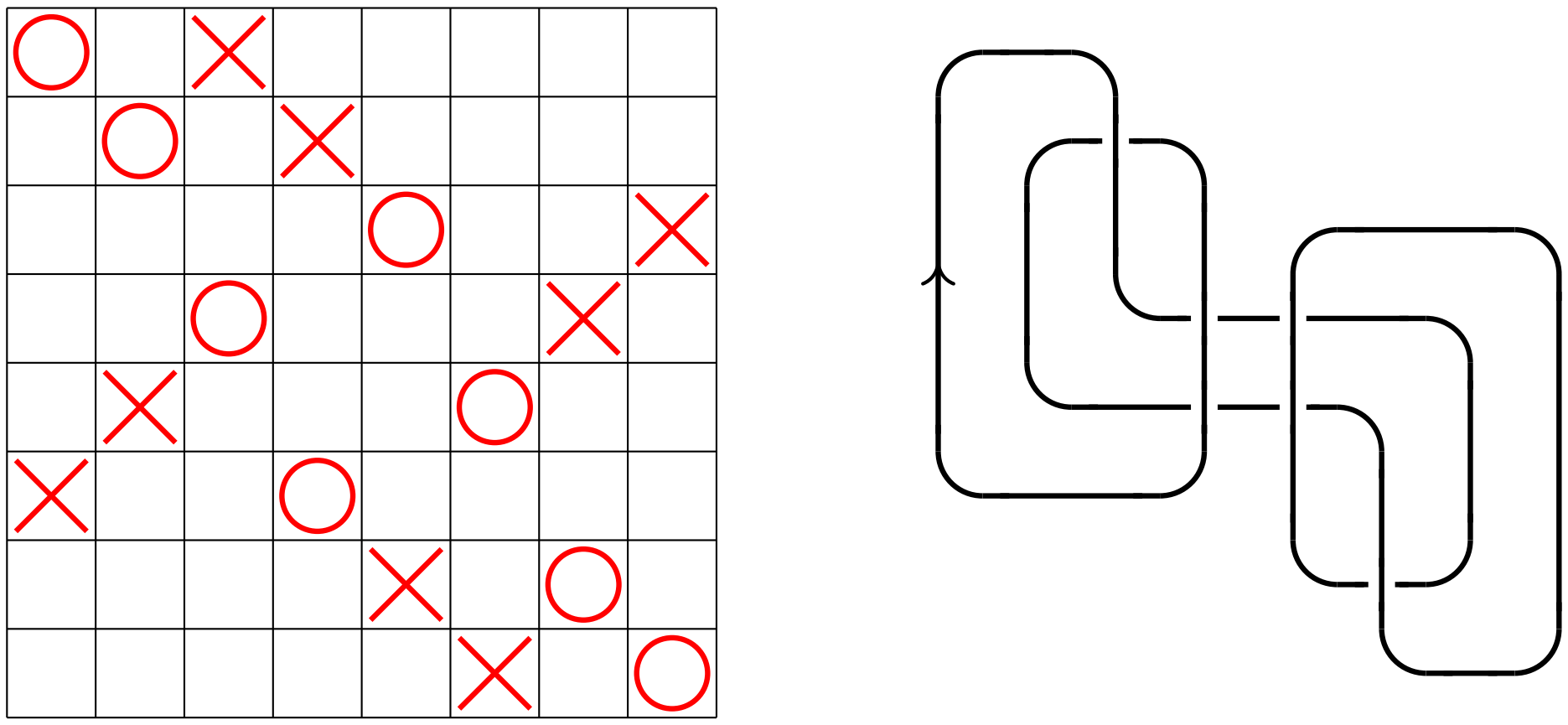}
	    \caption{A grid diagram with grid number $8$ and the oriented planar knot diagram that it determines.}
    	\label{fig:grid-diagram-example}
	\end{figure}
	
	\begin{definition}
	A \emph{grid diagram} with \emph{grid number} $n$ is an $n\times n$ square grid in which $n$ $X$s and $n$ $O$s have been placed in distinct squares, with each row and each column containing exactly one $X$ and one $O$.
	\end{definition}
	
	By connecting $X$ to $O$ in each column, connecting $O$ to $X$ in each row, declaring that vertical arcs pass over horizontal arcs, and smoothing corners in the obvious manner (see \Cref{fig:grid-diagram-example}), every grid diagram determines an oriented planar diagram for a smooth link. A fundamental result established in \cite{cromwell1995embedding} says that any two grid diagrams representing smoothly isotopic oriented knots or links are related by a sequence of modifications known as \emph{Cromwell moves}, playing for grid diagrams the role that Reidemeister moves play for planar diagrams of knots and links.
	
\begin{figure}[t]
	    \centering
	    \input{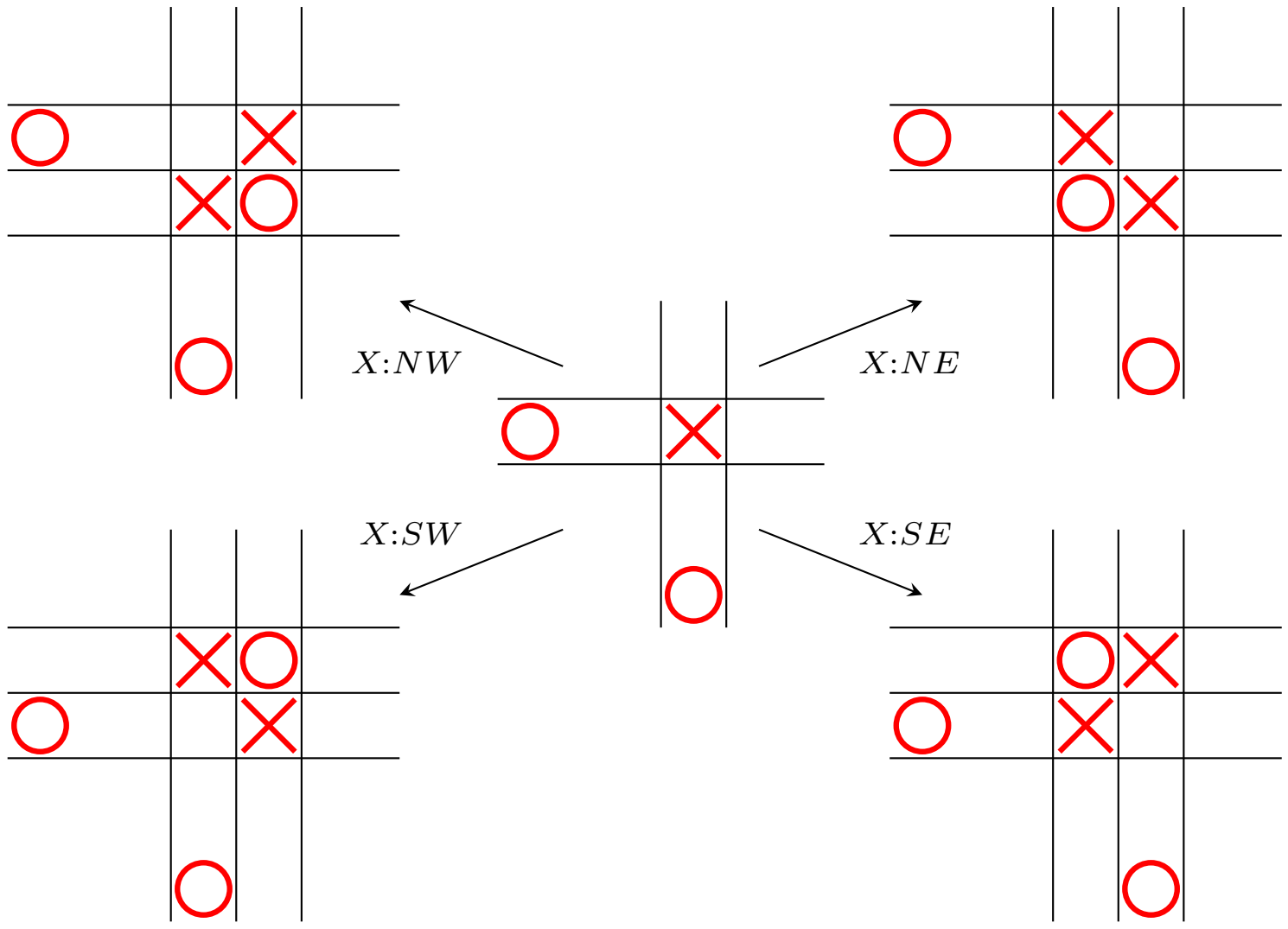}
	    \caption{Four stabilizations of a grid diagram at an $X$. This figure is adapted from \cite[Figure 1]{chongchitmate2013atlas}.}
    	    \label{fig:x-stabilizations}
	\end{figure}

	Cromwell moves, which we describe in greater detail in \cref{sec:background}, include:
	\begin{itemize}
		\item \emph{torus translations}, which move the topmost row of a grid diagram to the bottom;
		\item \emph{commutations}, which swap adjacent rows whose line segments are either disjoint or strictly nested when projected onto a single horizontal line;
		\item \emph{stabilizations}, which increase the grid number by 1, replacing a single $X$ with two $X$s and one $O$ in a manner depicted by \Cref{fig:x-stabilizations}.
	\end{itemize}
	The remaining Cromwell moves are obtained by rotating these by multiples of $90^\circ$; by swapping the roles of $X$ and $O$ in stabilization; and by inverting any Cromwell move. In \cite{dynnikov2006arc}, Dynnikov proved, among other things, the following monotonic simplification result for grid diagrams of the unknot.
	
	\begin{theorem}[{\cite[Theorem 1]{dynnikov2006arc}}]\label{thm:dynnikov-ms}
	For any grid diagram $D$ representing the unknot, there is a finite sequence of Cromwell moves which is nonincreasing on grid number and transforms $D$ into a grid diagram with grid number $2$.
	\end{theorem}
	
	An immediate consequence of this result is the following unknot detection algorithm: given any grid diagram $D$, there is a finite collection of grid diagrams obtained from $D$ by finite sequences of grid number-nonincreasing Cromwell moves, and an exhaustive search of this collection for a grid diagram with grid number $2$ determines whether or not $D$ represents the unknot.

	\subsection{Main results}
	Our first main result allows one to refine Dynnikov's unknot detection algorithm by considering only those Cromwell moves which preserve grid number.
	
	\begin{theorem}\label{thm:0-crossing}
	For any grid diagram $D$ representing the unknot, there is a finite sequence of torus translations and commutations which transforms $D$ into a $0$-crossing grid diagram.
	\end{theorem}
	
\begin{figure}[t]
	    \centering
	    \input{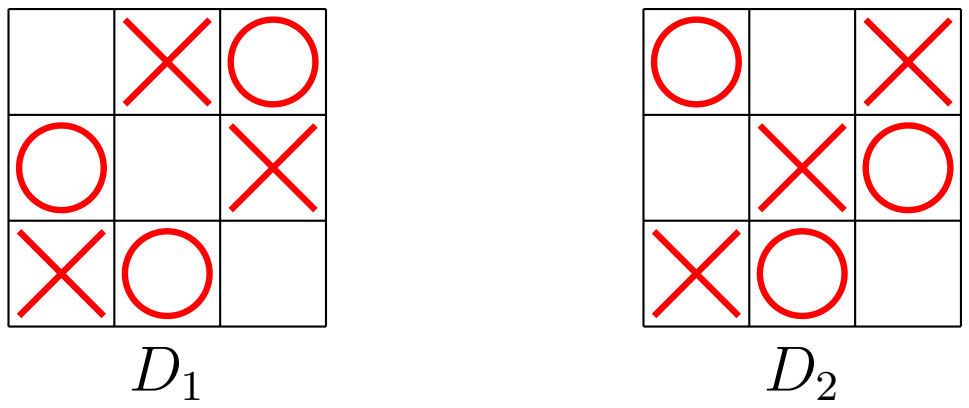}
	    \caption{Applying a degenerate commutation to the top two rows of $D_1$ produces $D_2$. These diagrams determine the same smooth isotopy type, but distinct Legendrian isotopy types.}
    	    \label{fig:degenerate-commutation}
	\end{figure}

	A result similar to \Cref{thm:0-crossing} was obtained by Hayashi-Yamada in \cite{hayashi2013unknotting}, with one crucial difference: while a commutation of adjacent rows (columns) requires these rows (columns) to be either disjoint or strictly nested when projected to a single line, Hayashi-Yamada allow \emph{degenerate commutations}, in which the projected rows (columns) may overlap at endpoints. See \Cref{fig:degenerate-commutation}.
	
	One reason for restricting attention to the usual (i.e., nondegenerate) Cromwell moves is that grid diagrams represent not only smooth isotopy classes, but also isotopy classes of braids; of Legendrian knots; and of transverse knots in $(S^3,\xi_{\mathrm{std}})$ (cf. \cite{ng2008grid}). While torus translations and commutations preserve all of these isotopy types (up to conjugation and exchange moves in the case of braids), their degenerate analogues need not do so. Indeed, by allowing degenerate Cromwell moves, Hayashi-Yamada obtain the following magnificent result:
	
	\begin{theorem}[{\cite[Theorem 1.2]{hayashi2013unknotting}}]\label{thm:hy}
	Any two grid diagrams of the unknot which have the same grid number are related by a sequence of possibly-degenerate torus translations and commutations.
	\end{theorem}
	
	It is immediately clear that a version of \Cref{thm:hy} which allows only nondegenerate torus translations and commutations cannot possibly be true. Indeed, for every grid number $n\geq 3$, there are grid diagrams representing Legendrian unknots in distinct Legendrian isotopy types; because torus translation and commutation preserve Legendrian isotopy type, these grid diagrams cannot be related by such moves. It is natural, therefore, to ask which grid diagrams \emph{are} related by sequences of torus translations and commutations. Let us denote by $\mathcal{G}_n$ the collection of grid diagrams of the unknot which have grid number $n\geq 2$ and consider the equivalence classes
	\[
	G_n := \mathcal{G}_n/\{\text{torus translations}, \text{commutations}\}.
	\]
	We will say that two grid diagrams in $\mathcal{G}_n$ are \emph{SPC-equivalent}\footnote{For ``size-preserving Cromwell."} if they represent the same element of $G_n$. Because grid diagrams determine oriented knots (with vertical arcs oriented from $X$ to $O$), we will say that two grid diagrams in $\mathcal{G}_n$ are \emph{unoriented SPC-equivalent} if they represent the same element of
	\[
	\overline{G}_n := \mathcal{G}_n/\{\text{torus translations}, \text{commutations}, \text{interchange of $X$s and $O$s}\}.
	\]
	Our next result identifies the number of SPC- and unoriented SPC-equivalence classes for each grid number.
	
	\begin{theorem}\label{thm:spc-classes}
	For $n\geq 2$, the sizes of $G_n$ and $\overline{G}_n$ are given by
	\[
	|G_n| = \dfrac{n(n^2-1)}{6}
	\quad\text{and}\quad
	|\overline{G}_n| =
		\begin{cases}
			\dfrac{n(n^2-1)}{12}, & \text{if \(n\) is odd},\\[8pt]
			\dfrac{n(n^2+2)}{12}, & \text{if \(n\) is even}.
		\end{cases}.
	\]
	\end{theorem}
	
	\begin{remark}
	For $n\geq 2$, the size of $G_n$ is recognized as the $(n-1)$th \emph{tetrahedral number}, while the sequence of sizes $(|\overline{G}_n|)_{n=2}^\infty$ appears as sequence \cite[A005993]{oeis}.
	\end{remark}
	
	\subsection{Applications}
	Beyond independent interest in SPC-equivalence, \Cref{thm:spc-classes} is practically useful for generating grid diagrams of the unknot. Consider the problem of enumerating all elements of $\mathcal{G}_n$, for some $n\geq 2$. \Cref{thm:dynnikov-ms,thm:hy,thm:spc-classes} provide the following solutions to this problem, respectively:
	\begin{enumerate}
		\item According to \Cref{thm:dynnikov-ms}, every element of $\mathcal{G}_n$ is obtained from a grid diagram of grid number $2$ by some sequence of torus translations, commutations, and $n-2$ stabilizations, with each stabilization increasing the grid number by $1$. We obtain $\mathcal{G}_n$ by considering all such sequences, deleting repeated diagrams.
		
		\item Using \Cref{thm:hy}, we may omit stabilizations, beginning with an arbitrary element of $\mathcal{G}_n$ and obtaining all others via sequences of possibly-degenerate torus translations and commutations.
		
		\item Finally, \Cref{thm:spc-classes} produces a similar algorithm. Beginning with a ``seed" representative for each equivalence class in $G_n$, we obtain all elements of $\mathcal{G}_n$ by applying sequences of torus translations and commutations to these seeds.
	\end{enumerate}
	While implementing the third algorithm requires generating multiple seed grid diagrams, the more general class of modifications used in the second algorithm leads to a greater number of repetitious grids. Moreover, the fact that the seeds in the third algorithm lie in distinct equivalence classes allows us to compare new grid diagrams against only those which were generated from the same seed, rather than against all grid diagrams generated to this point in the algorithm. For example, $\mathcal{G}_8$ includes $165,873,120$ distinct grid diagrams, and thus the final grid diagram generated by the second algorithm must be compared against at least $165,873,119$ previously-seen grids. However, these grid diagrams fall into 84 SPC-equivalence classes, ranging in size from $8$ grids to $10,775,808$ grids, meaning that the third algorithm compares each new grid to far fewer previously-seen grids.
	
	The algorithm just described for generating grid diagrams of the unknot has allowed us to study a more restrictive version of the monotonic simplification problem solved by Dynnikov in \Cref{thm:dynnikov-ms}. While \Cref{thm:dynnikov-ms} guarantees the existence of a simplifying sequence of Cromwell moves which is nonincreasing on grid number, another natural measure of complexity --- the crossing number --- is not considered. Inspired by other notions of \emph{hardness} for unknot diagrams (cf. \cite{applebaum2024unknotting,burton2024hard,gowers2011hard,hass2001reidemeister,hass2010quadratic,henrich2014unknotting,kauffman2011collapsing,petronio2016unknots}), we make the following definitions.
	
	\begin{definition}
	A grid diagram $D$ representing the unknot is said to be \emph{weakly Cromwell hard} if there is no sequence of Cromwell moves transforming $D$ into a grid diagram with grid number $2$ which is nonincreasing on both grid number and crossing number. We say that $D$ is \emph{strongly Cromwell hard} if it cannot be transformed into a grid diagram with grid number $2$ via a sequence of Cromwell moves which is nonincreasing on crossing number.
	\end{definition}
	
	Strong Cromwell hardness clearly implies weak Cromwell hardness. \jr{However, t}he latter notion is much more amenable to study, since a fixed grid diagram admits infinitely many simple simplifying sequences --- that is, simplifying sequences with no repeated states --- but only finitely many are nonincreasing on grid number. Indeed, in \Cref{sec:applications} we exhibit several grid diagrams for the unknot whose status as weakly Cromwell hard we have checked through exhaustive search. Using the algorithm described above to generate all unknot grid diagrams of grid number at most $8$, we have the following computational result.
	
	\begin{theorem}\label{thm:min-weakly-hard}
	The smallest weakly Cromwell hard grid diagrams for the unknot have grid number $8$.
	\end{theorem}
	
	In particular, we have identified a grid number- and crossing number-nonincreasing simplifying sequence of Cromwell moves for each grid diagram of the unknot with grid number at most $8$, with $128$ exceptions. For those $128$ exceptional grid diagrams, all of which have grid number $8$, we have verified by exhaustive search that no such sequence exists, making these the smallest (by grid number) weakly Cromwell hard grid diagrams.
	
	By restricting to a preferred collection of Legendrian, braid, or transverse isotopy types, \Cref{thm:min-weakly-hard} can be pushed further. Specifically, while the second algorithm above generates all of $\mathcal{G}_n$ from a single seed, SPC-equivalence allows us to generate only those elements of $\mathcal{G}_n$ representing our preferred Legendrian, braid, or transverse isotopy types, resulting in a more manageable computation. For instance, consider the following result.
	
	\begin{theorem}\label{thm:min-weakly-hard-max-tb}
	Every weakly Cromwell hard grid diagram for the Legendrian unknot with Thurston-Bennequin number $-1$ has grid number at least $9$.
	\end{theorem}
	
	The computations which verify \Cref{thm:min-weakly-hard,thm:min-weakly-hard-max-tb} are discussed in \Cref{sec:applications}, where we also exhibit a weakly Cromwell hard grid diagram for the Legendrian unknot with Thurston-Bennequin number $-1$ and grid number $14+4n$, for every $n\geq 1$. The \emph{Mathematica} files used to carry out these computations are available from the first author's website.
	
	\subsection{Organization}
	\Cref{sec:background} provides background on Cromwell moves and Legendrian invariants for grid diagrams. \Cref{thm:0-crossing} is established in \Cref{sec:proof-of-0-crossing-thm} and the proof of \Cref{thm:spc-classes} is split across \Cref{sec:staircases,sec:spc-equivalence-classes}. Finally, \Cref{sec:applications} applies our main results to the study of weakly Cromwell hard grid diagrams, establishing \Cref{thm:min-weakly-hard,thm:min-weakly-hard-max-tb}.

	\subsection*{Acknowledgments}
	This work was partially supported by NSF grant DMS-2532551. The authors are grateful to the Frost Fund in Cal Poly's Bailey College of Science and Mathematics for providing a supportive research environment.
	
	\subsection*{Declaration of AI use}
	The authors used ChatGPT to proofread this document (resulting in the identification of several minor errors), as well as to troubleshoot the \emph{Mathematica} code used to produce the computations in \Cref{sec:applications}. AI was not used for any other purpose.
\section{Background}\label{sec:background}
	
	Our primary object of interest, the grid diagram, was defined above to be an $n\times n$ square grid in which $n$ $X$s and $n$ $O$s have been placed, with each row and each column containing exactly one $X$ and one $O$. The present section has three goals: to establish our convention for presenting grid diagrams via their horizontal and vertical arcs; to provide a more careful treatment of stabilizations; and to define the Thurston-Bennequin and rotation numbers of a grid diagram.
	
	\subsection{Arc images}
	As explained in the introduction, our interest in grid diagrams is owed to their ability to represent knots and links in $\mathbb{R}^3$; we restrict attention to grid diagrams which represent knots. When it is appropriate to do so, we will therefore present grid diagrams not as grids filled with $X$s and $O$s, but as \emph{arc images}. Namely, each row and column of a grid diagram determines a single \emph{arc} in $\mathbb{R}^2$, which is a line segment connecting the $X$ in that row or column to the $O$ in that row or column, and the arc image of a grid diagram is the collection of these arcs. Notice that an arc image determines a grid diagram only up to interchanging $X$s and $O$s --- that is, only up to orientation. Despite this ambiguity, we will occasionally, when orientation is unimportant, refer to a grid diagram interchangeably with its arc image.
	
\begin{figure}
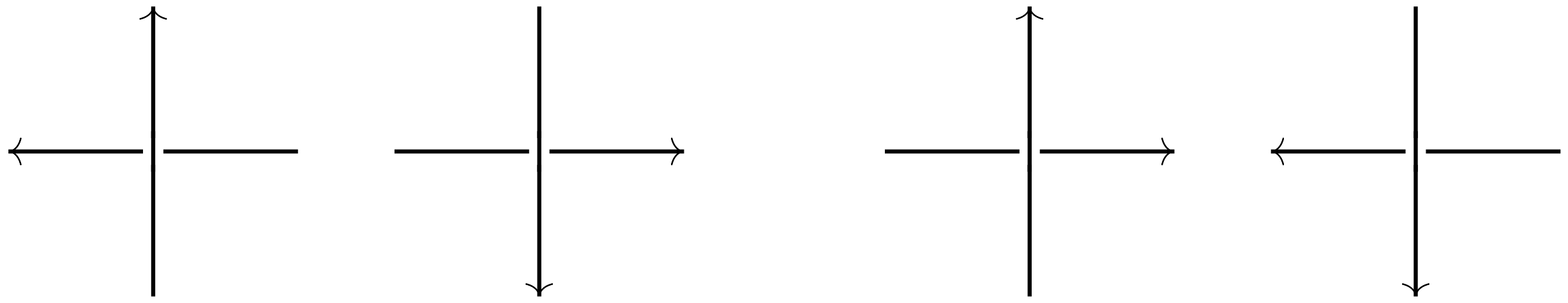

	\centering
		\begin{subfigure}{0.45\textwidth}
			\centering
			\input{positive-crossings.tex}
			\caption{positive crossings}
		\end{subfigure}\hfill
		\begin{subfigure}{0.45\textwidth}
			\centering
			\input{negative-crossings.tex}
			\caption{negative crossings}
		\end{subfigure}
		\caption{Positive and negative crossings}
		\label{fig:crossing-signs}
	\end{figure}

	Though the corners of an arc image are not rounded and the crossing data is not recorded, the arc image of a grid diagram is evidently equivalent to an unoriented knot diagram. In particular, the corners can be smoothed in a manner which is canonical up to planar isotopy, and our convention is that vertical strands always pass over horizontal strands at crossings. In this way, any value computed from a knot diagram can be computed from a grid diagram. For example, a grid diagram has a crossing wherever the interiors of two arcs in its arc image intersect. Using the convention that vertical arcs are oriented from $X$ to $O$ and horizontal arcs are oriented from $O$ to $X$, we may classify each crossing as either positive or negative in the usual manner, depicted in \Cref{fig:crossing-signs}. We may thus speak of the \emph{writhe} of a grid diagram, defined by
	\[
	\mathrm{wr}(D) := c_+(D) - c_-(D),
	\]
	where $c_+(D)$ and $c_-(D)$ give the number of positive and negative crossings of $D$, respectively.
	
	\subsection{Stabilizations}
\begin{figure}[t]
	    \centering
	    \input{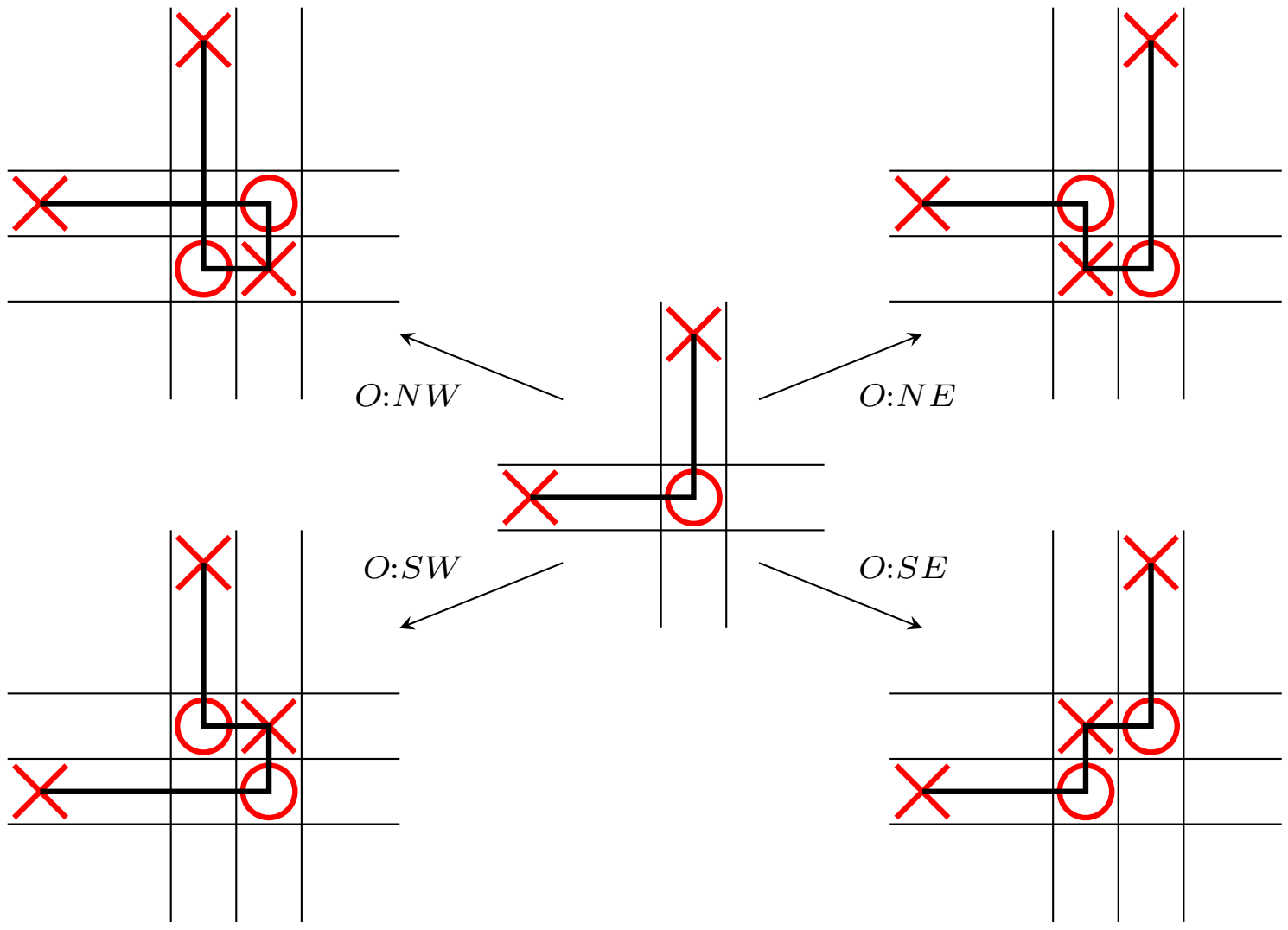}
	    \caption{The four stabilizations of a grid diagram at an $O$. The thick black arcs correspond to the arc image of the grid diagram.}
    	    \label{fig:o-stabilizations}
	\end{figure}

	Stabilizations are the only Cromwell moves to increase the grid number. A stabilization can be applied at any $X$ or $O$ of a grid diagram, replacing an $X$ with a $2\times 2$ square containing two $X$s and one $O$ or replacing an $O$ with a $2\times 2$ square containing two $O$s and one $X$. This gives eight types of stabilization, each denoted by the letter being replaced and the direction of the empty entry in the $2\times 2$ square which is introduced:
	\[
	\Xne, \Xnw, \Xsw, \Xse, \One, \Onw, \Osw, \Ose.
	\]
	\Cref{fig:x-stabilizations} shows the four possible stabilizations at an $X$ which is an upper-right corner of the grid; \Cref{fig:o-stabilizations} depicts the four possible stabilizations at an $O$ which is a lower-right corner of the grid. The Cromwell move which is inverse to a stabilization (and thus decreases the grid number by 1) is known as a \emph{destabilization}.
	
	At any corner, one of the four stabilizations introduces a crossing in the corresponding knot diagram, corresponding to a type I Reidemeister move, while the other three stabilizations preserve the planar isotopy type of the knot diagram. Of the planar-isotopy-type-preserving stabilizations, those with the appearance of ``concave dents" will be especially useful for our purposes.
	
	\begin{definition}
	A \emph{simple dent stabilization} is a stabilization of any of the following types:
	\begin{itemize}
		\item type $\Xne$ or $\One$ at an upper-right corner;
		\item type $\Xnw$ or $\Onw$ at an upper-left corner;
		\item type $\Xsw$ or $\Osw$ at a lower-left corner;
		\item type $\Xse$ or $\Ose$ at a lower-right corner.
	\end{itemize}
	\end{definition}
	
	For use in some of our later arguments, we also give a name to those stabilizations which introduce a crossing.
	
	\begin{definition}
		A \emph{twist stabilization} is a stabilization of any of the following types:
		\begin{itemize}
			\item type $\Xsw$ or $\Osw$ at an upper-right corner;
			\item type $\Xse$ or $\Ose$ at an upper-left corner;
			\item type $\Xne$ or $\One$ at a lower-left corner;
			\item type $\Xnw$ or $\Onw$ at a lower-right corner.
		\end{itemize}
	\end{definition}
	
	Stabilizations provide a convenient way to construct grid diagrams of the unknot with large grid number, leading to the following definition.
	
	\begin{definition}
	We say that a grid diagram $D$ is \emph{purely-stabilized} if $D$ can be obtained from a grid diagram of grid number $2$ via a sequence of stabilizations, excluding other Cromwell moves. If each of these stabilizations can be taken to be a simple dent stabilization, then we say that $D$ is a \emph{purely dent-stabilized grid diagram}.
	\end{definition}
	
	\begin{remark}
	Observe that a purely-stabilized grid diagram necessarily represents the unknot.
	\end{remark}
	
	Finally, we establish notation for a particular construction of purely dent-stabilized grid diagrams. Throughout, we number rows top-to-bottom and columns left-to-right, as in a matrix.
	
	\begin{definition}
	Given $a,b,c,d\in \mathbb{Z}_{\geq 0}$, we define the purely dent-stabilized grid diagrams $D_+(a,b,c,d)$ and $D_-(a,b,c,d)$ recursively as follows:
	\begin{itemize}
		\item $D_+(0,0,0,0)$ is the unique grid diagram of grid number $2$ with an $X$ in row $1$, column $1$; and $D_-(0,0,0,0)$ is the grid obtained from $D_+(0,0,0,0)$ by interchanging the roles of $X$ and $O$.
		
		\item For $a>0$, $D_\pm(a,0,0,0)$ is obtained from $D_\pm(a-1,0,0,0)$ by applying a $NW$ stabilization at row $1$, column $a$.
		
		\item For $b>0$ and $a\geq 0$, $D_\pm(a,b,0,0)$ is obtained from $D_\pm(a,b-1,0,0)$ by applying a $NE$ stabilization at row $b$, column $a+b+1$.
		
		\item For $c>0$ and $a,b\geq 0$, $D_\pm(a,b,c,0)$ is obtained from $D_\pm(a,b,c-1,0)$ by applying a $SW$ stabilization at row $a+b+c+1$, column $c$.
		
		\item For $d>0$ and $a,b,c\geq 0$, $D_\pm(a,b,c,d)$ is obtained from $D_\pm(a,b,c,d-1)$ by applying a $SE$ stabilization at row $a+b+c+d+1$, column $a+b+c+2$.
	\end{itemize}
	The shared arc image of $D_+(a,b,c,d)$ and $D_-(a,b,c,d)$ is denoted by $D(a,b,c,d)$.
	\end{definition}
	
	In plain language, $D_\pm(a,b,c,d)$ is obtained by applying $a$ simple dent stabilizations in the upper-left corner of a grid diagram of grid number $2$; then applying $b$ simple dent stabilizations in the upper-right corner; then $c$ simple dent stabilizations in the lower-left corner; and finally $d$ simple dent stabilizations in the lower-right corner. Notice that $D_\pm(a,b,c,d)$ has grid number $a+b+c+d+2$.
	
	\begin{definition}
	The grid diagrams $D_\pm(a,b,c,d)$, for any choice of $a,b,c,d\geq 0$, are said to be in \emph{standard form}.
	\end{definition}
	
	The proofs of our main theorems amount to verifying that every grid diagram of the unknot is SPC-equivalent to a grid diagram in standard form, and then counting the SPC-equivalence classes of grid diagrams in standard form. For examples of standard form grid diagrams, see \Cref{fig:pds-to-std,fig:reversal-lemma}.
	
	\subsection{Classical Legendrian invariants}
	As mentioned above, grid diagrams can be used to represent Legendrian knots in the standard contact 3-manifold $(\mathbb{R}^3,\xi_{\mathrm{std}})$. For the purposes of the present paper, no knowledge of contact geometry is necessary; instead, we may think of Legendrian knots as combinatorial objects. In particular, \cite[Proposition 2]{ng2008grid} demonstrates that the set of Legendrian isotopy classes of Legendrian knots is in bijective correspondence with the set
	\[
	\mathcal{L} := \{\text{grid diagrams}\}/\sim_{\mathcal{L}},
	\]
	where $\sim_{\mathcal{L}}$ denotes equivalence up to SPC moves, $\Xne$ (de)stabilization, and $\Xsw$ (de)stabilization. Given a grid diagram $D$, we will refer to the equivalence class in $\mathcal{L}$ represented by $D$ as the Legendrian isotopy type of $D$.
	
	\begin{figure}
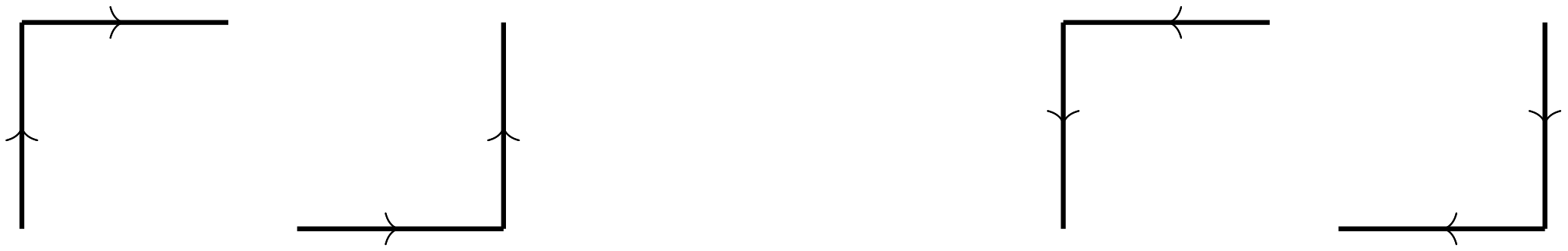

	\centering
		\begin{subfigure}{0.45\textwidth}
			\centering
			\input{up-cusps.tex}
			\caption{up cusps}
		\end{subfigure}\hfill
		\begin{subfigure}{0.45\textwidth}
			\centering
			\input{down-cusps.tex}
			\caption{down cusps}
		\end{subfigure}
		\caption{Up and down cusps}
		\label{fig:cusps}
	\end{figure}
	
	We will adapt to this convention two integer-valued invariants of Legendrian isotopy classes which are known as the \emph{classical invariants}. Given a grid diagram $D$, we refer to the upper-left and lower-right corners of $D$ as the \emph{left cusps} and \emph{right cusps} of $D$, respectively. A cusp of $D$ is said to be an \emph{up cusp} if the vertical arc meeting that cusp is oriented from bottom to top and is said to be a \emph{down cusp} if this arc is oriented from top to bottom. See \Cref{fig:cusps}. With this vocabulary, we may define the \emph{Thurston-Bennequin invariant} of $D$, written $\tb(D)$, to be
	\[
	\tb(D) := \mathrm{wr}(D) - \tfrac{1}{2}\mathrm{cusp}(D),
	\]
	where $\mathrm{cusp}(D)$ is the total number of cusps of $D$, and the \emph{rotation number} of $D$, written $\rot(D)$, to be
	\[
	\rot(D) := \tfrac{1}{2}(\mathrm{down}(D)-\mathrm{up}(D)),
	\]
	where $\mathrm{down}(D)$ and $\mathrm{up}(D)$ record the number of down cusps and up cusps of $D$, respectively. One may verify that the Thurston-Bennequin invariant and rotation number are invariants of Legendrian isotopy classes by checking that both are preserved under SPC moves as well as $\Xne$ and $\Xsw$ stabilization.
\section{Crossingless grid diagrams}\label{sec:proof-of-0-crossing-thm}
	Throughout this section, we treat all grid diagrams as unoriented, and thus depict only their arc images. The purpose of this section is to prove \Cref{thm:0-crossing}, which we deduce from the following stronger result:
	
	\begin{theorem}\label{thm:strongtheorem1}
	Every grid diagram of the unknot is SPC-equivalent to a purely dent-stabilized diagram in standard form.
	\end{theorem}
	
	\Cref{thm:strongtheorem1} will be deduced from a sequence of lemmas which show that standard-form grid diagrams provide representations of the SPC-equivalence classes of progressively larger collections of grid diagrams. The sequence of implications is summarized as follows:
	\begin{align*}
	\left\{\begin{matrix} \text{grid diagrams}\\ \text{of the unknot}\end{matrix}\right\}
	&\xrightarrow[\text{\Cref{lem:stabForm}}]{\text{SPC-equivalent to}}
	\left\{\begin{matrix} \text{purely-stabilized}\\ \text{grid diagrams}\end{matrix}\right\}\\
	&\xrightarrow[\text{\Cref{lem:onestabpuredent}}]{\text{SPC-equivalent to}}
	\left\{\begin{matrix} \text{purely dent-stabilized}\\ \text{grid diagrams}\end{matrix}\right\}\\
	&\xrightarrow[\text{\Cref{lem:pds-to-std}}]{\text{SPC-equivalent to}}
	\left\{\begin{matrix} \text{standard form}\\ \text{grid diagrams}\end{matrix}\right\}
	\end{align*}

	As observed by Dynnikov (cf. \cite[Remark 12]{dynnikov2006arc}), in any sequence of Cromwell moves which is nonincreasing on grid number, it is simple to delay any destabilizations until the end of the sequence. Because we will use this technique throughout the arguments which follow, we provide a careful justification of this observation, which we encode via the following \emph{sequence lifting lemma}.
	
	\begin{lemma}\label{lem:sequenceLifting}
	For any sequence of SPC moves $D_1 \longrightarrow D_2 \longrightarrow \cdots \longrightarrow D_m$ and any grid diagram $G_1$ which is obtained from $D_1$ via stabilization, there exists a sequence of SPC moves $G_1 \longrightarrow G_2 \longrightarrow \cdots \longrightarrow G_M$ and a sequence of integers $1=i_1 < i_2 < \cdots <i_m = M$ such that, for $1\leq k\leq m$, $G_{i_k}$ is a stabilization of $D_k$. Moreover, if each SPC move $D_i\longrightarrow D_{i+1}$, $1\leq i < m$, is a commutation, then
	\begin{enumerate}
		\item each SPC move $G_j\longrightarrow G_{j+1}$ may be taken to be a commutation;
		\item if, additionally, $G_1$ is obtained from $D_1$ via simple dent stabilization, then $G_M$ may be obtained from $D_m$ via simple dent stabilization.
	\end{enumerate}
	\end{lemma}
	
	\begin{proof}
	It suffices to prove this result in the case where $m=2$. Given a sequence of SPC moves $D_1\longrightarrow D_2\longrightarrow \cdots \longrightarrow D_m$, with $m>2$, and a stabilization $G_1$ of $D_1$, we may then apply the $m=2$ case to the subsequence $D_1\longrightarrow D_2$, producing $G_1\longrightarrow G_2\longrightarrow \cdots \longrightarrow G_{i_2}$, with $G_{i_2}$ a stabilization of $D_2$. Iterating this process for the sequences $D_k\longrightarrow D_{k+1}$ and grids $G_{i_{k}}$, for $2\leq k< m$, constructs the desired sequence $G_1\longrightarrow G_2\longrightarrow \cdots \longrightarrow G_M$.
	
	Now suppose that $D_1\longrightarrow D_2$ is an SPC move, where $D_1$ has grid number $n$, and that $G_1$ is obtained from $D_1$ by a stabilization in row $r_0$ and column $c_0$, for some $1\leq r_0,c_0\leq n$. If $D_1\longrightarrow D_2$ commutes rows $i$ and $i+1$, for some $1\leq i < r_0-1$ or commutes columns $j$ and $j+1$, for some $1\leq j < c_0-1$, then we may produce $G_2$ from $G_1$ by commuting the same pair of rows or columns. On the other hand, if $D_1\longrightarrow D_2$ commutes rows $i$ and $i+1$, for some $r_0+1\leq i < n$, then we may define the move $G_1\longrightarrow G_2$ to be the commutation which swaps rows $i+1$ and $i+2$. Likewise, if $D_1\longrightarrow D_2$ swaps columns $j$ and $j+1$, for some $c_0+1 \leq j < n$, then the move $G_1\longrightarrow G_2$ can be defined as the commutation of columns $j+1$ and $j+2$. In all of these cases, we are lifting a commutation of $D_1$ which does not interact with the stabilization of $D_1$ that produced $G_1$.
	
	Now if $D_1\longrightarrow D_2$ swaps row $r_0$ with either row $r_0-1$ or $r_0+1$, we will define a sequence $G_1\longrightarrow G_2\longrightarrow G_3$ of two commutations, with $G_3$ a stabilization of $D_2$. Consider the case where $D_1\longrightarrow D_2$ swaps rows $r_0$ and $r_0-1$. Because this commutation is possible, the arcs contained in rows $r_0$ and $r_0-1$ must be either properly nested or disjoint when projected to a parallel axis, and this property is preserved by stabilization. Namely,
	\begin{itemize}
		\item if the arc of row $r_0$ is nested within that of row $r_0-1$ in $D_1$, then the arcs of rows $r_0$ and $r_0+1$ are nested within that of row $r_0-1$ in $G_1$;
		\item if the arc of row $r_0$ is disjoint from that of row $r_0-1$ in $D_1$, then the arcs of rows $r_0$ and $r_0+1$ are disjoint from that of row $r_0-1$ in $G_1$;
		\item and if the arc of row $r_0-1$ is nested within that of row $r_0$ in $D_1$, then the arc of row $r_0-1$ is nested within the arc of either row $r_0$ or $r_0+1$ in $G_1$, and is disjoint from the arc of the other row.
	\end{itemize}
	For the final case above, see \Cref{fig:stabilized-rows}. In any case, we may define the sequence $G_1\longrightarrow G_2\longrightarrow G_3$ to be the commutation of rows $r_0-1$ and $r_0$ followed by the commutation of rows $r_0$ and $r_0+1$. With this definition, $G_3$ is the result of stabilizing $D_2$ at row $r_0-1$ and column $c_0$. In the case where $D_1\longrightarrow D_2$ swaps rows $r_0$ and $r_0+1$, we may define $G_1\longrightarrow G_2\longrightarrow G_3$ to swap row $r_0+2$ with row $r_0+1$ and then swap row $r_0+1$ with row $r_0$; justifying that these commutations are possible follows as in the case described here.
	
	\begin{figure}
	\centering
	\input{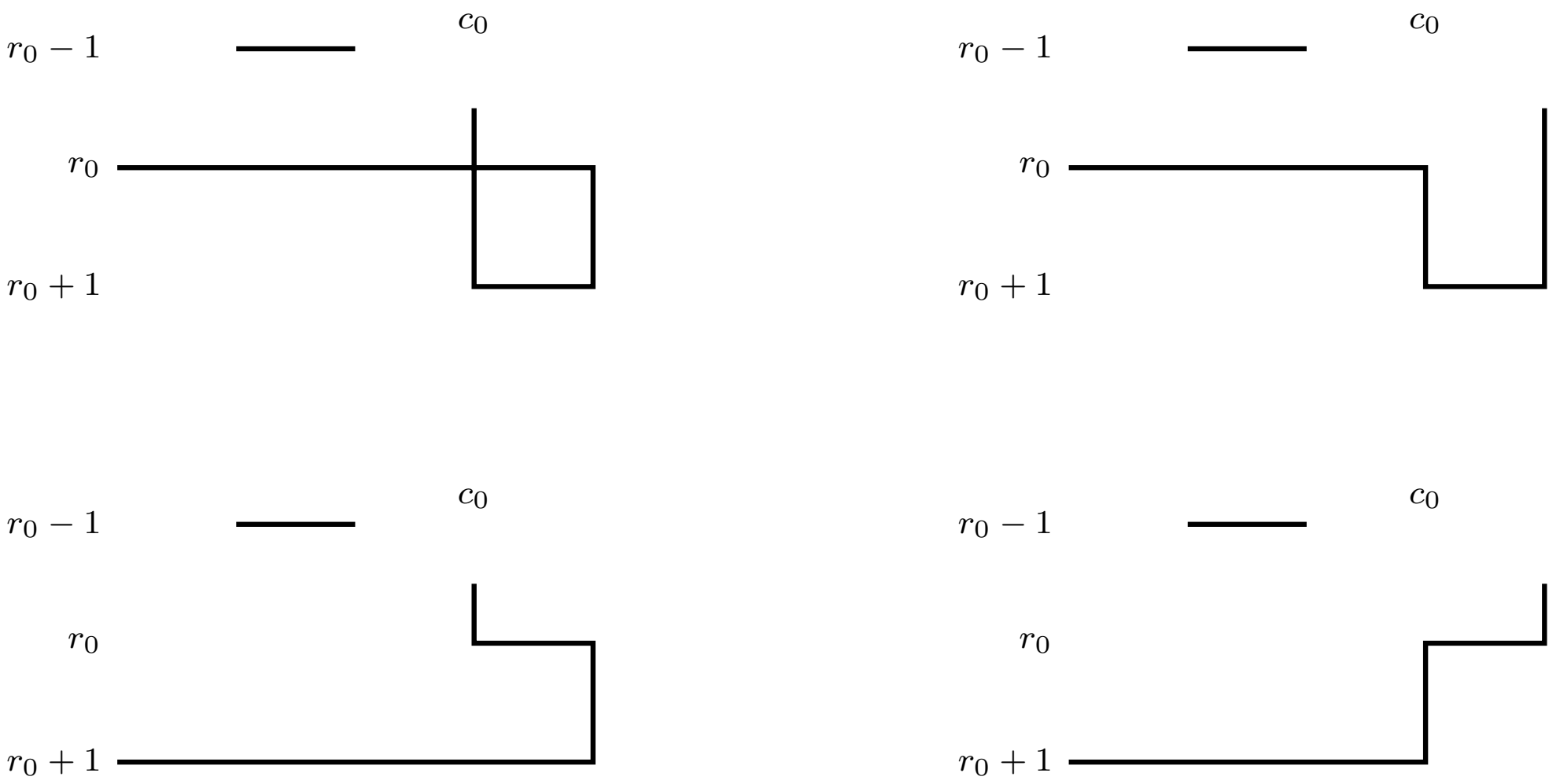}
	\caption{If the arc of row $r_0-1$ is nested within the arc of row $r_0$ before stabilizing in row $r_0$, then the arc of row $r_0-1$ will be disjoint from the ``new" arc and nested within the ``old" arc after stabilizing.}
	\label{fig:stabilized-rows}
	\end{figure}
	
	The same argument applies if $D_1\longrightarrow D_2$ swaps column $c_0$ with either column $c_0-1$ or $c_0+1$, and thus it remains to consider the case where $D_1\longrightarrow D_2$ is a torus translation. Specifically, let us suppose that $D_1\longrightarrow D_2$ is the torus translation which moves the arc in row $1$ to row $n$; the other three torus translations will follow from an analogous argument. Now if $r_0>1$, then we may simply define $G_1\longrightarrow G_2$ to be the torus translation which moves the arc in row $1$ to row $n+1$; we then see that $G_2$ is obtained from $D_2$ by stabilizing in row $r_0-1$ and column $c_0$. On the other hand, if $r_0=1$, then we define the sequence $G_1\longrightarrow G_2\longrightarrow G_3$ to be \emph{two} top-to-bottom torus translations --- that is, $G_1\longrightarrow G_2$ moves the arc in row $1$ to row $n+1$ and $G_2\longrightarrow G_3$ does the same. The grid diagram $G_3$ is then obtained from $D_2$ by a stabilization in row $n$ and column $c_0$. In particular, this stabilization has the same type as that which was performed on $D_1$ in row $1$ and column $c_0$ to produce $G_1$.
	
	Altogether, we see that each SPC move $D_k\longrightarrow D_{k+1}$, for $1\leq k < m$, lifts to a sequence of SPC moves $G_{i_k}\longrightarrow G_{i_k+1}\longrightarrow \cdots \longrightarrow G_{i_{k+1}}$, with $G_{i_{k+1}}$ a stabilization of $D_{k+1}$. Inductively applying this to the moves of $D_1\longrightarrow D_2\longrightarrow \cdots \longrightarrow D_m$ produces the desired sequence $G_1\longrightarrow G_2\longrightarrow \cdots \longrightarrow G_M$.
	\end{proof}

	Using the sequence lifting lemma, we now establish that, in proving \Cref{thm:strongtheorem1}, it suffices to show that a grid diagram is SPC-equivalent to a purely dent-stabilized grid diagram.
	
\begin{figure}
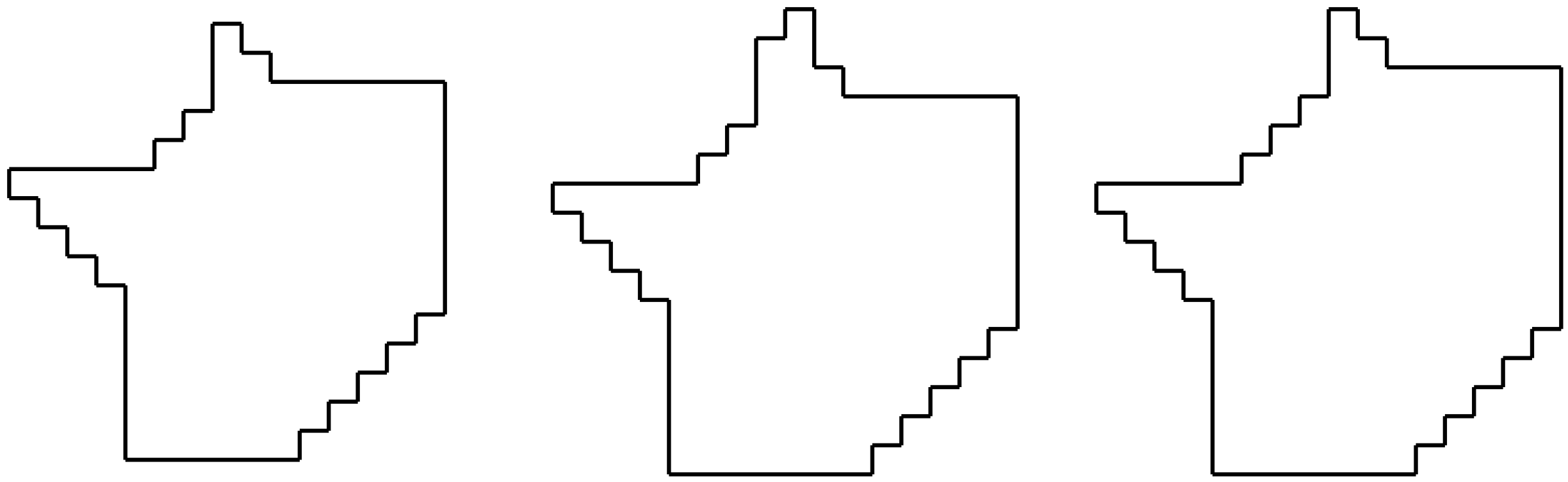

	\centering
		\begin{subfigure}{0.32\textwidth}
			\centering
			\input{pds-to-std-1.tex}
		\end{subfigure}\hfill
		\begin{subfigure}{0.32\textwidth}
			\centering
			\input{pds-to-std-2.tex}
		\end{subfigure}\hfill
		\begin{subfigure}{0.32\textwidth}
			\centering
			\input{pds-to-std-3.tex}
		\end{subfigure}
		\caption{The first panel shows $D(3,2,4,5)$. Applying a simple dent stabilization at row $1$, column $8$ produces a purely dent-stabilized grid diagram which is not in standard form, seen in the second panel. A sequence of $2$ commutations produces $D(4,2,4,5)$, seen in the third panel.}
		\label{fig:pds-to-std}
	\end{figure}

	\begin{lemma}\label{lem:pds-to-std}
	Any purely dent-stabilized grid diagram is commutation-equivalent to a purely dent-stabilized grid diagram in standard form.
	\end{lemma}
	\begin{proof}
	Let $G$ be a purely dent-stabilized grid diagram of grid number $n$. Our argument proceeds by induction on $n$. There are precisely two grid diagrams of grid number $2$ --- $D_+(0,0,0,0)$ and $D_-(0,0,0,0)$ --- and these provide a base case for our induction. Next, we suppose that $n>2$, and that any purely dent-stabilized grid diagram of grid number less than $n$ is commutation-equivalent to a purely dent-stabilized grid diagram in standard form.
	
	Because $G$ is purely dent-stabilized, there is a purely dent-stabilized grid diagram $D$ of grid number $n-1$ such that $G$ is obtained from $D$ by a simple dent stabilization. According to our inductive hypothesis, $D$ is commutation-equivalent to some $D_\pm(a,b,c,d)$, with $a+b+c+d=n-3$, and \Cref{lem:sequenceLifting} then tells us that, up to commutation-equivalence, $G$ is obtained from $D_\pm(a,b,c,d)$ via simple dent stabilization. There are four positions in $D_\pm(a,b,c,d)$ where applying a simple dent stabilization might fail to produce a grid diagram in standard form. They are:
	\begin{itemize}
		\item row $1$, column $a+c+1$;
		\item row $a+b+1$, column $1$;
		\item \jr{row $n$}, column $c+1$;
		\item row $b+1$, \jr{column $n$}.
	\end{itemize}
	In the first and third cases, simple dent stabilization introduces a horizontal arc of length $1$ which may be commuted down to row $b+2$ (in the first case) or up to row $a+b+c+2$ (in the third case) until we have $D_\pm(a+1,b,c,d)$ (in the first case) or $D_\pm(a,b,c+1,d)$ (in the third case). See \Cref{fig:pds-to-std}. The second and fourth cases similarly produce grid diagrams which are commutation-equivalent to $D_\pm(a+1,b,c,d)$ and $D_\pm(a,b+1,c,d)$, respectively. We conclude that $G$ is commutation-equivalent to a grid diagram in standard form, as desired.
	\end{proof}
	
	Next, we use the sequence lifting lemma to reduce the proof of \Cref{thm:strongtheorem1} to the case of purely-stabilized grid diagrams.
	
	\begin{lemma}\label{lem:stabForm}
		Every grid diagram of the unknot is SPC-equivalent to a purely-stabilized grid diagram.
	\end{lemma}
	
	\begin{proof}
	Let $D$ be a grid diagram of grid number $n$ representing the unknot. From \Cref{thm:dynnikov-ms} we obtain a sequence
	\begin{equation}\label{eq:dynnikov-ms}
	D=D_0 \longrightarrow D_1 \longrightarrow \cdots \longrightarrow D_{\ell+(n-2)}=U_2
	\end{equation}
	of Cromwell moves which are nonincreasing on grid number that transforms $D$ into a $2\times 2$ grid diagram. We claim that we may produce from this sequence a new sequence
	\begin{equation}\label{eq:delayed-destabilizations}
	D=D'_0 \longrightarrow D'_1 \longrightarrow \cdots \longrightarrow D'_{\ell'} \longrightarrow \cdots \longrightarrow D'_{\ell'+(n-2)}=U_2
	\end{equation}
	of Cromwell moves with the property that $D'_{k-1}\longrightarrow D'_{k}$ is an SPC move for $1\leq k\leq \ell'$ and $D'_k\to D'_{k+1}$ is a destabilization for $\ell' \leq k< \ell'+(n-2)$. That is, we may delay the destabilizations in \Cref{eq:dynnikov-ms} until the end of the sequence. To justify this claim, it suffices to show that if $D_{k-1}\longrightarrow D_k\longrightarrow D_{k+1}$ consists of a destabilization followed by an SPC move, then there is another sequence
	\[
	D_{k-1}=G_0\longrightarrow G_1\longrightarrow \cdots \longrightarrow G_{M}\longrightarrow G_{M+1} = D_{k+1}
	\]
	of Cromwell moves, where $G_{j-1}\longrightarrow G_j$ is an SPC move, for $1\leq j\leq M$. Repeated application of this sub-claim to \Cref{eq:dynnikov-ms} then produces \Cref{eq:delayed-destabilizations}. But the sub-claim is an immediate consequence of applying \Cref{lem:sequenceLifting} to the (rather short) sequence of SPC moves given by $D_k\longrightarrow D_{k+1}$ and the stabilization $D_{k-1}$ of $D_k$. Namely, \Cref{lem:sequenceLifting} provides the sequence of SPC moves $G_0\longrightarrow G_1\longrightarrow\cdots\longrightarrow G_M$, with the guarantee that $G_M$ is a stabilization of $D_{k+1}$.
	\end{proof}
	
	In light of \Cref{lem:pds-to-std,lem:stabForm}, we can establish \Cref{thm:strongtheorem1} by showing that every purely-stabilized grid diagram is SPC-equivalent to a purely dent-stabilized grid diagram. For this, the sequence lifting lemma affords an inductive strategy wherein it will suffice to show that the result of applying any stabilization to a purely dent-stabilized grid diagram is SPC-equivalent to a purely dent-stabilized diagram.

	\begin{lemma}\label{lem:onestabpuredent}
	Applying a single stabilization to a purely dent-stabilized grid diagram produces a grid diagram which is SPC-equivalent to a purely dent-stabilized grid diagram.
	\end{lemma}
	\begin{proof}
	We will denote by $D$ a purely dent-stabilized grid diagram and by $G$ a grid diagram obtained from $D$ via a single stabilization, and show that $G$ is SPC-equivalent to a purely dent-stabilized grid diagram. Notice that if $G$ is obtained from $D$ via a dent stabilization, then this conclusion is immediate.
	
	\begin{figure}
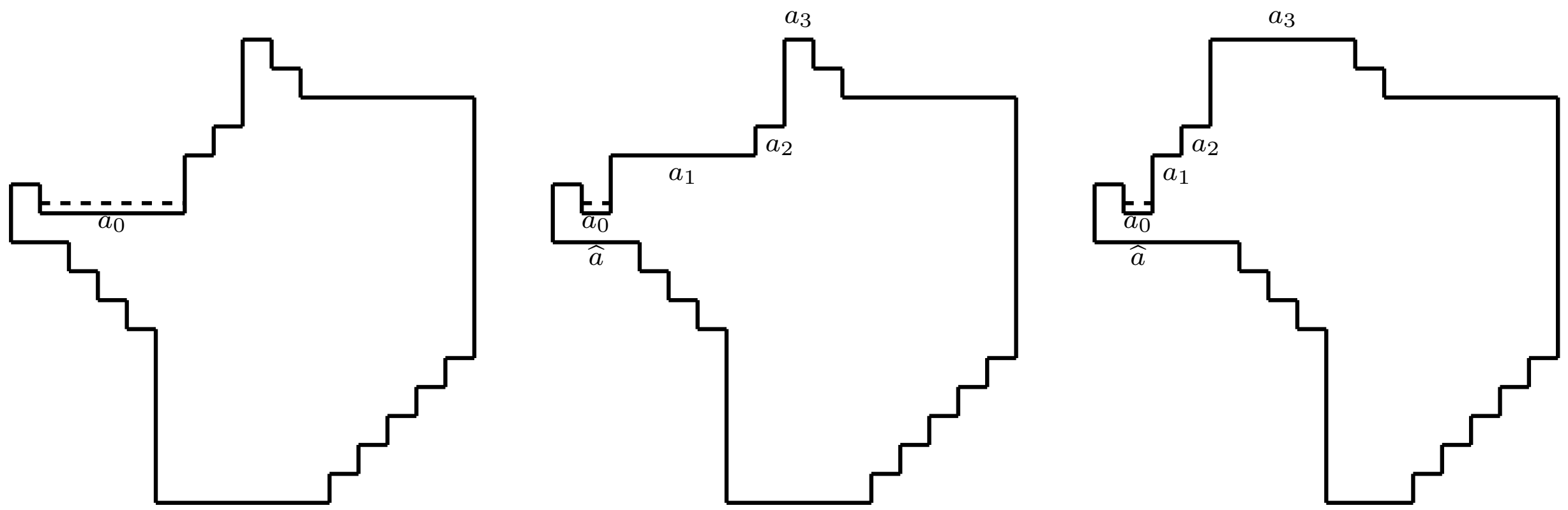

	\centering
		\begin{subfigure}{0.32\textwidth}
			\centering
			\input{goofy-crossingless-1.tex}
			\caption{The arc $a_0$ which witnesses non-convexity has $\ell(a_0)>1$.}
		\end{subfigure}\hfill
		\begin{subfigure}{0.32\textwidth}
			\centering
			\input{goofy-crossingless-2.tex}
			\caption{Commutations ensure that $\ell(a_0)=1$.}
		\end{subfigure}\hfill
		\begin{subfigure}{0.32\textwidth}
			\centering
			\input{goofy-crossingless-3.tex}
			\caption{Further commutations put $a_0,a_1,\ldots,a_{m-1}$ ``within" $\widehat{a}$.}
			\label{subfig:goofy-crossingless-3}
		\end{subfigure}
		\caption{A crossingless stabilization which is not a dent (performed here in the leftmost column) introduces a concavity. In this figure, we apply commutations to ensure that the sequence of arcs $a_0,a_1,\ldots,a_{m-1}$ may be simultaneously commuted past the arc $\widehat{a}$.}
		\label{fig:goofy-crossingless}
	\end{figure}

	Next, we claim that if $G$ is obtained from $D$ via a crossingless stabilization which is not a dent stabilization, then $G$ is SPC-equivalent to a grid diagram $G'$ which is obtained from some purely dent-stabilized $D'$ via a dent stabilization. To justify the claim, we first invoke \Cref{lem:pds-to-std} to observe that $D$ is commutation-equivalent to $D(a,b,c,d)$, for some $a,b,c,d\geq 0$. After applying a corresponding sequence of commutations to $G$ (via the sequence lifting lemma), we may assume that $G$ is obtained from $D(a,b,c,d)$ via a single crossingless stabilization. If $G$ is \emph{axis-convex} --- meaning that every horizontal or vertical line segment with endpoints in the (closed) region $R(G)$ in $\mathbb{R}^2$ bounded by $G$ is contained in $R(G)$ --- then $G$ is evidently commutation-equivalent to a purely dent-stabilized grid diagram. If $G$ fails to be axis-convex, then there is a unique arc $a_0$ in $G$ with the property that a small horizontal or vertical translation of $a_0$ produces an arc whose endpoints lie on the arc image of $G$ and whose interior points lie outside $R(G)$. Moreover, we may apply commutations to ensure that $a_0$ has length $1$, while its endpoints lie on arcs $b_1$ and $c_1$ of length $1$ and $2$, respectively. Finally, there is a unique arc $\widehat{a}$ in $G$ which is parallel to $a_0$ and has the property that, when projected to a parallel axis, the endpoints of $a_0$ are properly nested within those of $\widehat{a}$.
	
	\begin{figure}
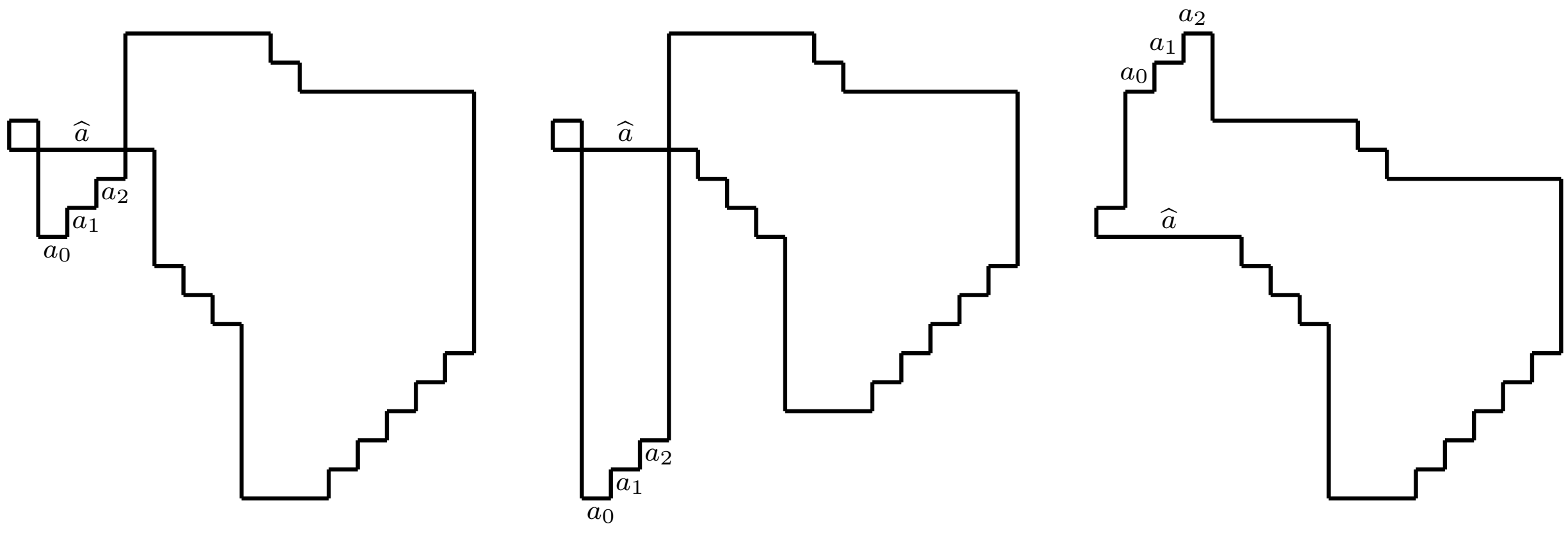

	\centering
		\begin{subfigure}{0.32\textwidth}
			\centering
			\input{past-a-hat.tex}
		\end{subfigure}\hfill
		\begin{subfigure}{0.32\textwidth}
			\centering
			\input{to-extreme-row.tex}
		\end{subfigure}\hfill
		\begin{subfigure}{0.32\textwidth}
			\centering
			\input{after-tt.tex}
		\end{subfigure}
		\caption{The first panel is obtained by commuting the arcs $a_0,a_1,a_2$ in \Cref{subfig:goofy-crossingless-3} past $\widehat{a}$; the second by further commuting these arcs into an extreme position; and the third by applying three bottom-to-top torus translations.}
		\label{fig:crossingless-spc}
	\end{figure}
	
	Now consider the maximal sequence of parallel arcs $a_0,a_1,\ldots,a_m$ with the following features: $a_1$ has an endpoint on $c_1$; for $1\leq i \leq m$, the arcs $a_{i-1}$ and $a_i$ share one endpoint when projected to a parallel axis; and for $1\leq i\leq m$, the arc connecting an endpoint of $a_{i-1}$ to an endpoint of $a_i$ is oriented in the same manner as is $c_1$. We may, if necessary, perform commutations to ensure that the endpoints of the arcs $a_0,a_1,\ldots,a_{m-1}$ are contained in the interior of $\widehat{a}$ when projected to a parallel axis. See \Cref{fig:goofy-crossingless} for a case where such commutations are necessary. Following this normalization, each of the arcs $a_0,a_1,\ldots,a_{m-1}$ may be commuted past $\widehat{a}$ as in the first panel of \Cref{fig:crossingless-spc}. Moreover, we may apply further commutations to ensure that, once past $\widehat{a}$, $a_0$ lies in an extreme row or column and the arc connecting $a_{i-1}$ to $a_i$ has length $1$, for $1\leq i\leq m-1$, as in the second panel of \Cref{fig:crossingless-spc}. Using the same names for the arcs following these commutations, we next apply $m$ torus translations to bring the arcs $a_0,a_1,\ldots,a_{m-1}$ across the grid diagram and return $\widehat{a}$ to its original position, as seen in the final panel of \Cref{fig:crossingless-spc}. Following these torus translations, we have a zero-crossing grid which is axis-convex, and thus is commutation-equivalent to a purely dent-stabilized grid diagram. We see that performing a crossingless stabilization on a purely dent-stabilized grid diagram produces, up to SPC moves, a purely dent-stabilized grid diagram.
		
	Finally, let us consider the case where $G$ is obtained from $D$ via a twist stabilization. Again, using \Cref{lem:sequenceLifting,lem:pds-to-std}, we may assume that $D$ has the form $D(a,b,c,d)$, for some $a,b,c,d\geq 0$; we also assume that $D$ has grid number $n$. We will show that $G$ is SPC-equivalent to a grid diagram $G'$ which is obtained from some purely dent-stabilized grid diagram $D'$ via a crossingless stabilization. This will complete our proof, as the case of crossingless stabilizations has been covered above.
	
	\begin{figure}
	\centering
	\input{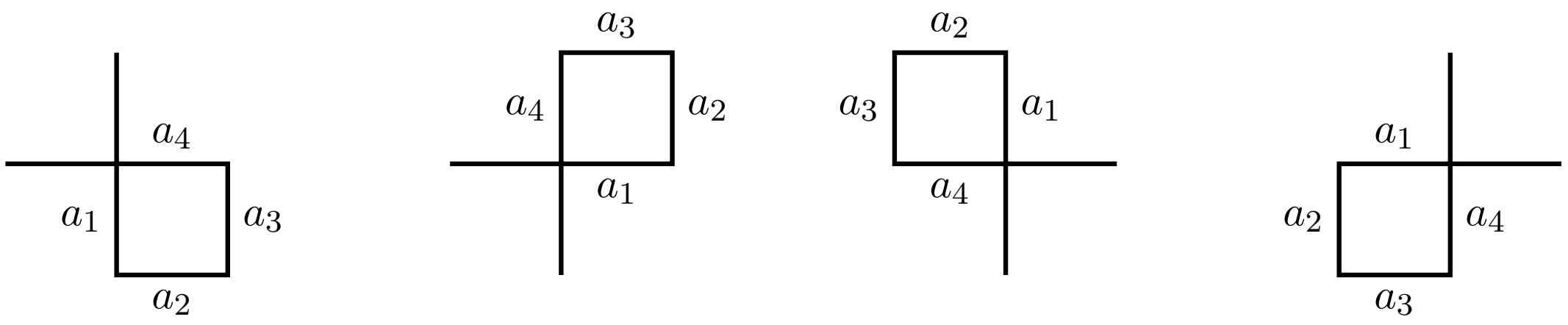}
	\caption{From left to right, twist stabilizations of type $NW$, $SW$, $SE$, and $NE$, respectively. In the first two cases, the \emph{guiding arc} is $a_2$; in the latter two cases, the guiding arc is $a_3$.}
	\label{fig:twist-stab-arcs}
	\end{figure}

	A twist stabilization naturally produces a $1\times 1$ region in $\mathbb{R}^2$ bounded by arcs of $G$. Beginning from the crossing and proceeding in a counterclockwise manner, let us label these arcs $a_1,a_2,a_3,a_4$. Notice that $a_2$ and $a_3$ each have length $1$, while $a_1$ and $a_4$ do not. Let us denote by $\gamma$ the \emph{guiding arc} of the stabilization, which we define to be $a_2$ if we have performed a stabilization of type $NW$ or $SW$ and define to be $a_3$ if we have performed a stabilization of type $NE$ or $SE$. See \Cref{fig:twist-stab-arcs}. This arc will guide a sequence of SPC moves which transforms our twist stabilization into a stabilization of another type. As a final preparation for the process which produces this sequence, define the pair $(r,c)$ so that $r$ is the index of the topmost row containing one of $a_1,a_2,a_3$, or $a_4$; and $c$ is the index of the leftmost column containing one of $a_1,a_2,a_3$, or $a_4$.
	
	\begin{figure}
	\centering
	\input{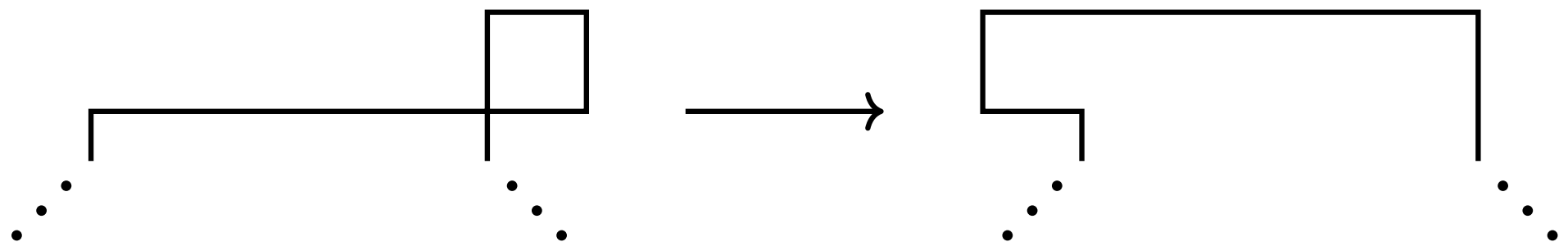}
	\caption{By repeatedly commuting its guiding arc to the right and performing a single right-to-left torus translation, a twist stabilization of type $SW$ applied in the top row of a diagram is transformed into a crossingless stabilization.}
	\label{fig:sw-twist-top-row}
	\end{figure}
	
	Now suppose our twist stabilization is of type $SW$ and consider the following sequence of SPC moves: we commute the guiding arc $\gamma$, which begins in column $c+1$, to the right until $\gamma$ is in column $n+1$; then apply a right-to-left torus translation; and then continue commuting $\gamma$ to the right as far as possible. If $r=1$, then the resulting grid diagram is obtained from a purely dent-stabilized grid diagram via a crossingless stabilization, as depicted in \Cref{fig:sw-twist-top-row}. On the other hand, if $r>1$, then we may commute $\gamma$ to the right until it lies in column $c-1$; notice that $c\geq 2$ for any twist stabilization of type $SW$. In this case, our new grid diagram is obtained from a purely dent-stabilized grid via a twist stabilization of type $NE$ and the value of $c$ has been decremented by $1$.
	
	If our twist stabilization is of type $NE$, we perform an analogous sequence of SPC moves, in this case commuting the guiding arc $\gamma$, which begins in row $r+1$, down to row $n+1$; then bringing $\gamma$ back to the top with a torus translation; and finally commuting $\gamma$ down as far as possible. Similar to the previous case, this process will produce a grid diagram obtained from a purely dent-stabilized grid diagram via crossingless stabilization if $c=1$. If $c>1$, then we will perform $n-1$ SPC moves in total, bringing $\gamma$ to row $r-1\geq 1$ and ending with a grid diagram which is obtained from a purely dent-stabilized grid diagram via a twist stabilization of type $SW$. Similar to the case of a $SW$ twist stabilization, this process decrements the value of $r$ by $1$.
		
	\begin{figure}
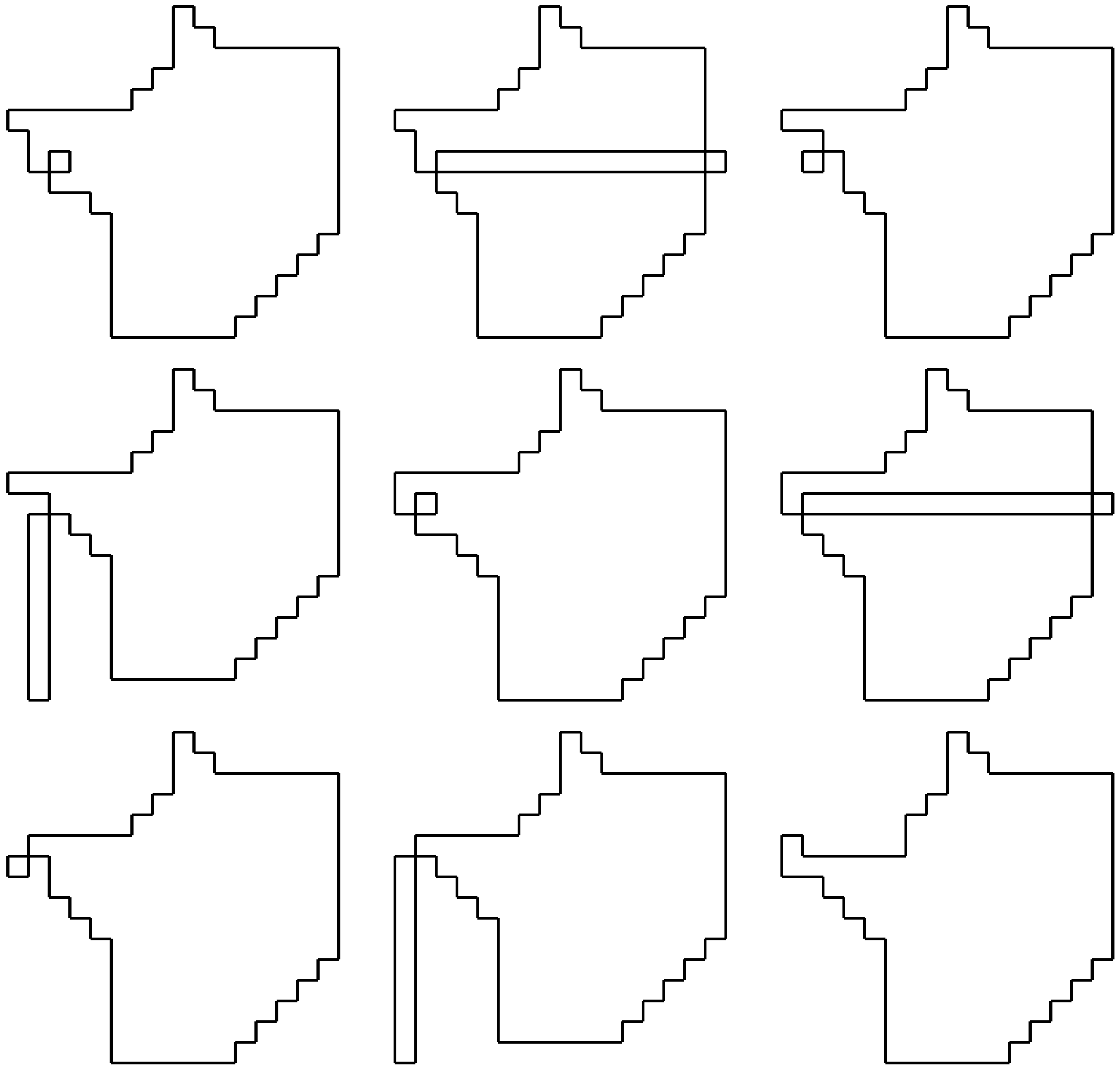

	\centering
		\begin{subfigure}{0.32\textwidth}
			\centering
			\input{twisty-movie-1.tex}
		\end{subfigure}\hfill
		\begin{subfigure}{0.32\textwidth}
			\centering
			\input{twisty-movie-2.tex}
		\end{subfigure}\hfill
		\begin{subfigure}{0.32\textwidth}
			\centering
			\input{twisty-movie-3.tex}
		\end{subfigure}
		
		\vspace{1em}
		
		\begin{subfigure}{0.32\textwidth}
			\centering
			\input{twisty-movie-4.tex}
		\end{subfigure}\hfill
		\begin{subfigure}{0.32\textwidth}
			\centering
			\input{twisty-movie-5.tex}
		\end{subfigure}\hfill
		\begin{subfigure}{0.32\textwidth}
			\centering
			\input{twisty-movie-6.tex}
		\end{subfigure}
		
		\vspace{1em}
		
		\begin{subfigure}{0.32\textwidth}
			\centering
			\input{twisty-movie-7.tex}
		\end{subfigure}\hfill
		\begin{subfigure}{0.32\textwidth}
			\centering
			\input{twisty-movie-8.tex}
		\end{subfigure}\hfill
		\begin{subfigure}{0.32\textwidth}
			\centering
			\input{twisty-movie-9.tex}
		\end{subfigure}
		\caption{In the top-left, the result of applying a twist stabilization of type $SW$ to a purely dent-stabilized grid diagram. The top row shows the first $n-1$ SPC moves, producing a twist stabilization of type $NE$. Another sequence of $n-1$ SPC moves produces a twist stabilization of type $SW$, closer to the leftmost column of the grid diagram. The process is repeated until we obtain the bottom-right grid diagram, which results from applying a crossingless stabilization to a purely dent-stabilized grid diagram.}
		\label{fig:twisty-movie}
	\end{figure}
	
	In this way, we see that applying a twist stabilization of type $SW$ or $NE$ to a purely dent-stabilized grid diagram is equivalent, up to SPC moves, to performing a crossingless stabilization: namely, by alternating between the two cases just described, we construct a sequence of SPC moves which decrements the values of $r$ and $c$ until one of them is $1$; the next step in the process provides an SPC equivalence to a grid obtained via crossingless stabilization on a purely dent-stabilized grid diagram. See \Cref{fig:twisty-movie}. An analogous algorithm exists for twist stabilizations of type $SE$ or $NW$, and thus, combined with our earlier treatment of crossingless stabilizations, we see that the grid diagram which results from applying a twist stabilization of any type to a purely dent-stabilized grid diagram is SPC-equivalent to a purely dent-stabilized grid diagram. This concludes the proof.
	\end{proof}

We are now prepared to use \Cref{lem:sequenceLifting,lem:pds-to-std,lem:stabForm,lem:onestabpuredent} to establish \Cref{thm:strongtheorem1}.

	\begin{proof}[Proof of \Cref{thm:strongtheorem1}]
	According to \Cref{lem:stabForm}, every grid diagram of the unknot is SPC-equivalent to a purely-stabilized grid diagram, and \Cref{lem:pds-to-std} tells us that purely dent-stabilized grid diagrams are SPC-equivalent to grid diagrams in standard form. Thus we need only prove that purely-stabilized grid diagrams are SPC-equivalent to purely dent-stabilized grid diagrams. We do so by induction on the grid number $n$. A purely-stabilized grid diagram of grid number $3$ is obtained from the unique $2\times 2$ grid diagram via one stabilization and thus, by \Cref{lem:onestabpuredent}, is SPC-equivalent to a purely dent-stabilized grid diagram. Now suppose the result holds for all purely stabilized grid diagrams of grid number $n\geq 3$. If $G$ is a grid diagram of grid number $n+1$, then we identify a purely-stabilized grid diagram $D$ of grid number $n$ of which $G$ is a stabilization. Our inductive hypothesis ensures that $D$ is transformed by a sequence of SPC-moves into a purely dent-stabilized grid diagram, and \Cref{lem:sequenceLifting} allows us to lift this sequence to a sequence of SPC-moves applied to $G$. This shows that $G$ is SPC-equivalent to a grid diagram obtained from a purely dent-stabilized grid diagram by a single stabilization; finally, another application of \Cref{lem:onestabpuredent} produces an SPC-equivalence to a purely dent-stabilized grid diagram.
	\end{proof}
	
	In addition to implying \Cref{thm:0-crossing}, \Cref{thm:strongtheorem1} will be useful in our proof of \Cref{thm:spc-classes}, which we undertake in the next section.
\section{Staircases and inverted staircases}\label{sec:staircases}
	Our proof of \Cref{thm:spc-classes} will use SPC moves to systematically convert crossingless grid diagrams into purely dent-stabilized grid diagrams in standard form. These grid diagrams frequently contain certain sub-grids which we refer to as \emph{staircases} and \emph{inverted staircases}, and in this section we show that staircases and inverted staircases may be interchanged for one another via SPC moves. Throughout, we work with grid diagrams, rather than their arc images.
	
\begin{figure}
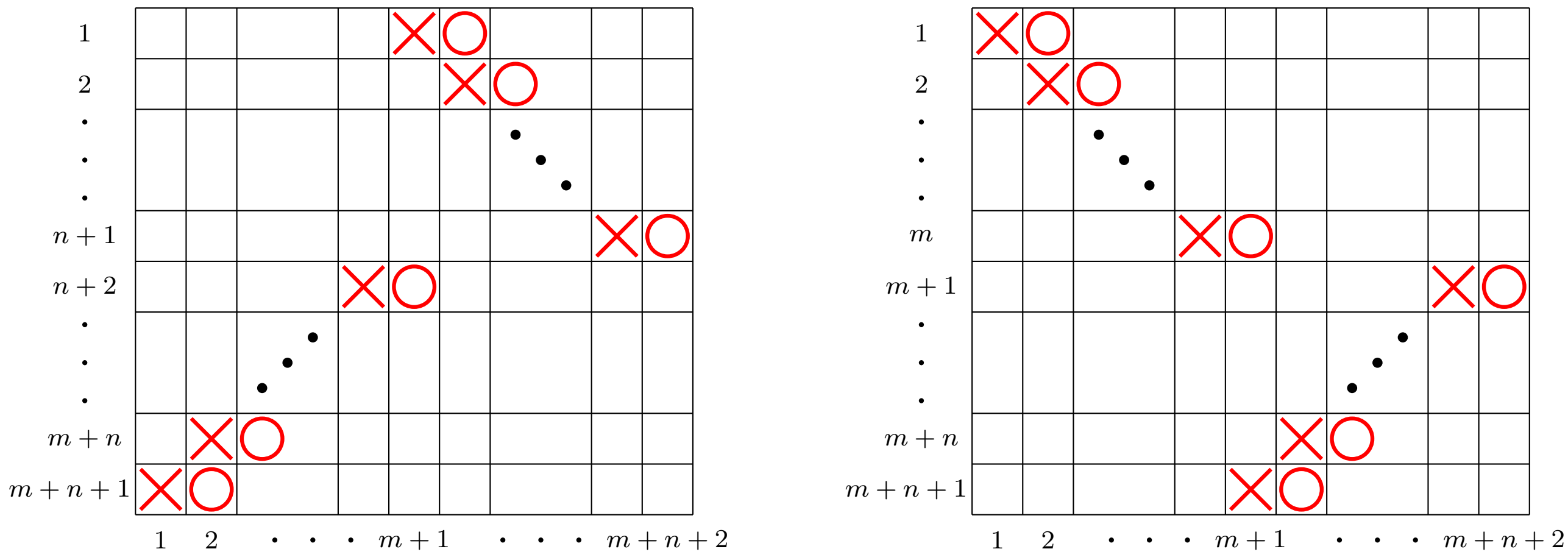

	\centering
	\begin{subfigure}{0.45\textwidth}
		\centering\tikzsetnextfilename{positive-staircase}
		\input{positive-staircase}
		\caption{The positive staircase $S^+_{m,n}$}
		\label{subfig:staircase}
	\end{subfigure}\hfill
	\begin{subfigure}{0.45\textwidth}
		\centering\tikzsetnextfilename{positive-inverted-staircase}
		\input{positive-inverted-staircase}
		\caption{The positive inverted staircase $Z^+_{m,n}$}
		\label{subfig:inverted-staircase}
	\end{subfigure}
	\caption{The positive staircase and positive inverted staircase naturally appear as sub-grids of a grid diagram.}
	\label{fig:staircases}
	\end{figure}

	Given integers $m,n\geq 0$, we will denote by $S^+_{m,n}$ and $Z^+_{m,n}$ the \emph{positive staircase} and \emph{positive inverted staircase} grids, respectively, defined as follows. Both are grids of size $(m+n+1)\times(m+n+2)$ (though neither is a grid diagram). For each integer $1\leq k\leq n$, row $k+1$ of $S^+_{m,n}$ has an $X$ in column $m+k+1$ and an $O$ in column $m+k+2$, while row $m+k$ of $Z^+_{m,n}$ has an $X$ in column $m+n+2-k$ and an $O$ in column $m+n+3-k$. For each integer $1\leq k\leq m$, row $n+k+1$ of $S^+_{m,n}$ has an $X$ in column $m+1-k$ and an $O$ in column $m+2-k$, while row $k$ of $Z^+_{m,n}$ has an $X$ in column $k$ and an $O$ in column $k+1$. Finally, row $1$ of $S^+_{m,n}$ has an $X$ in column $m+1$ and an $O$ in column $m+2$, while row $m+n+1$ of $Z^+_{m,n}$ has an $X$ in column $m+1$ and an $O$ in column $m+2$. See \Cref{fig:staircases}. If we write $S^+_{m,n}(r,c)$ for the entry of $S^+_{m,n}$ in row $r$ and column $c$, with a similar definition for $Z^+_{m,n}(r,c)$, then this definition can be summarized by
	\[
	S^+_{m,n}(r,c) = \left\{\begin{matrix}
	X, & \text{if~} 1\leq r\leq n+1 \text{~and~} c-r=m\\
	O, & \text{if~} 1\leq r\leq n+1 \text{~and~} c-r=m+1\\
	X, & \text{if~} n+2\leq r\leq m+n+1 \text{~and~} c+r=m+n+2\\
	O, & \text{if~} n+2\leq r\leq m+n+1 \text{~and~} c+r=m+n+3
	\end{matrix}\right.
	\]
	and
	\[
	Z^+_{m,n}(r,c) = \left\{\begin{matrix}
	X, & \text{if~} 1\leq r\leq m \text{~and~} c-r=0\\
	O, & \text{if~} 1\leq r\leq m \text{~and~} c-r=1\\
	X, & \text{if~} m+1\leq r\leq m+n+1 \text{~and~} c+r=2m+n+2\\
	O, & \text{if~} m+1\leq r\leq m+n+1 \text{~and~} c+r=2m+n+3
	\end{matrix}\right.,
	\]
	where the relevant grid has no entry in position $(r,c)$ if $(r,c)$ matches none of the listed conditions. We may similarly define the \emph{negative staircase} $S^-_{m,n}$ and \emph{negative inverted staircase} $Z^-_{m,n}$, obtained from $S^+_{m,n}$ and $Z^+_{m,n}$, respectively, by swapping the roles of $X$ and $O$.
	
	As defined, each row of $S^\pm_{m,n}$ and $Z^\pm_{m,n}$ contains an arc of length $1$. We may widen either $S^\pm_{m,n}$ or $Z^\pm_{m,n}$ by inserting empty columns between any two of the columns defined above. Concretely, if $I=(i_1,i_2,\ldots,i_{m+n+1})$ is a multi-index of non-negative integers, then we will denote by $S^\pm_{m,n,I}$ the grid obtained from $S^\pm_{m,n}$ by inserting $i_j$ empty columns between columns $j$ and $j+1$ in $S^\pm_{m,n}$, for $1\leq j\leq m+n+1$, and similarly define $Z^\pm_{m,n,I}$. Notice that $S^\pm_{m,n,I}$ and $Z^\pm_{m,n,I}$ are grids of size $(m+n+1)\times(m+n+2+|I|)$, where $|I|$ is the sum of the entries of $I$. At last, we have a result which allows us to invert staircases.

	\begin{proposition}\label{prop:inverted-staircase}
	Let $D$ be a grid diagram which contains $Z^\pm_{m,n,I}$ as a sub-grid, for some $m,n\geq 0$ and some multi-index $I=(i_1,i_2,\ldots,i_{m+n+1})$ of non-negative integers. Let $D'$ be the grid diagram obtained from $D$ by replacing this sub-grid with $S^\pm_{n,m,I}$. Then $D$ and $D'$ are SPC-equivalent.
	\end{proposition}
	\begin{proof}
	Let us first establish the case where $i_j=0$ for $1\leq j\leq m+n+1$. Notice that, by construction, $Z^\pm_{m,0}=S^\pm_{0,m}$ and $Z^\pm_{0,n}=S^\pm_{n,0}$, for all $m,n\geq 0$. In particular, both $Z^+_{0,0}$ and $S^+_{0,0}$ are grids of size $1\times 2$ containing an $X$ in column $1$ and an $O$ in column $2$, while the letters are reversed for $Z^-_{0,0}$ and $S^-_{0,0}$. These provide the base cases for an argument by induction on $\min\{m,n\}$. For our inductive hypothesis, fix an integer $k\geq 0$ and suppose that whenever $D$ is a grid diagram which contains $Z^\pm_{m',n'}$ as a sub-grid, for some $m',n'\geq 0$ with $\min\{m',n'\}\leq k$, $D$ is SPC-equivalent, via vertical commutations and bottom-to-top torus translations, to the grid diagram obtained by replacing this sub-grid with $S^\pm_{n',m'}$.
	
	\begin{figure}
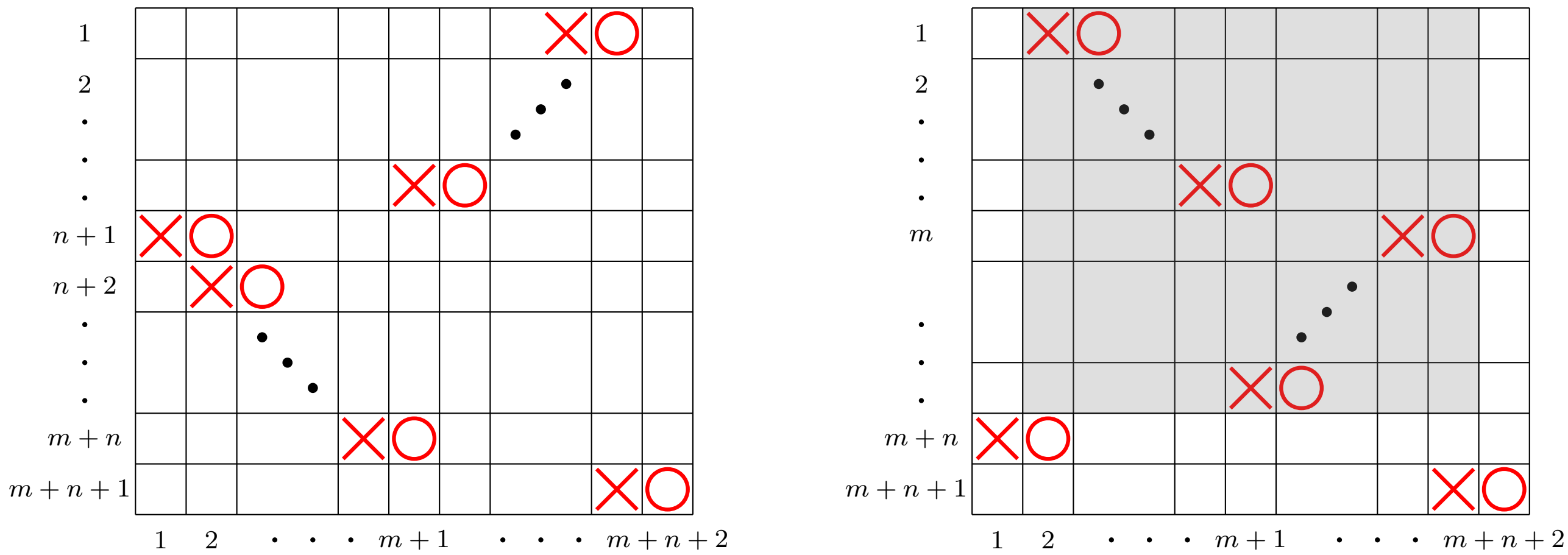

	\centering
	\begin{subfigure}{0.45\textwidth}
		\centering\tikzsetnextfilename{positive-inverted-staircase-step-1}
		\input{positive-inverted-staircase-step-1}
		\caption{After commuting the ``upward steps"}
		\label{subfig:inverted-step-1}
	\end{subfigure}\hfill
	\begin{subfigure}{0.45\textwidth}
		\centering\tikzsetnextfilename{positive-inverted-staircase-step-2}
		\input{positive-inverted-staircase-step-2}
		\caption{After commuting the ``downward steps"}
		\label{subfig:inverted-step-2}
	\end{subfigure}
	\caption{Our proof of \Cref{prop:inverted-staircase} begins by using commutations and torus translations to move all but one ``upward step" of $Z^\pm_{m,n}$ from the bottom of the grid to the top, then does the same to all but one ``downward step." After these moves, we recognize $Z^\pm_{m-1,n-1}$ as a sub-grid, shaded on the right.}
	\label{fig:staircase-spc1}
	\end{figure}
	
	Now consider a grid diagram $D$ which contains $Z^\pm_{m,n}$ as a sub-grid, for some $m,n\geq 0$ with $\min\{m,n\}=k+1$. We denote this sub-grid by $G$ and, as we apply SPC moves, continue to use $D$ to denote the (modified) grid diagram and $G$ to denote the sub-grid occupying these rows and columns. Because the arc in row $m+n+1$ of $G$ has length $1$, we may apply commutations which bring this arc down to the bottommost row of $D$. We may then apply a similar set of commutations which bring the arc originally in row $m+n$ of $G$ down to the second-from-bottom row of $D$. We continue in this fashion until the $n$ arcs which began in rows $m+2,m+3,\ldots,m+n+1$ of $G$ are in the bottom $n$ rows of $D$. Next, we apply $n$ bottom-to-top torus translations, moving these arcs to the top of $D$, and finally apply commutations, if necessary, to bring the $n$ arcs which began in rows $m+2,m+3,\ldots,m+n+1$ of $G$ into rows $1,2,\ldots,n$ of $G$. See \Cref{subfig:inverted-step-1}.
	
	Following the sequence of SPC moves just described, rows $n+1,n+2,\ldots,m+n$ of $G$ each contain a length-$1$ arc. In a fashion similar to the sequence of commutations, torus translations, and commutations described above, we now commute the $m-1$ arcs in rows $n+2,n+3,\ldots,m+n$ down to the bottom of $D$, apply $m-1$ bottom-to-top torus translations, and then commute these arcs back down into rows $1,2,\ldots,m-1$ of $G$, as depicted in \Cref{subfig:inverted-step-2}.

	\begin{figure}
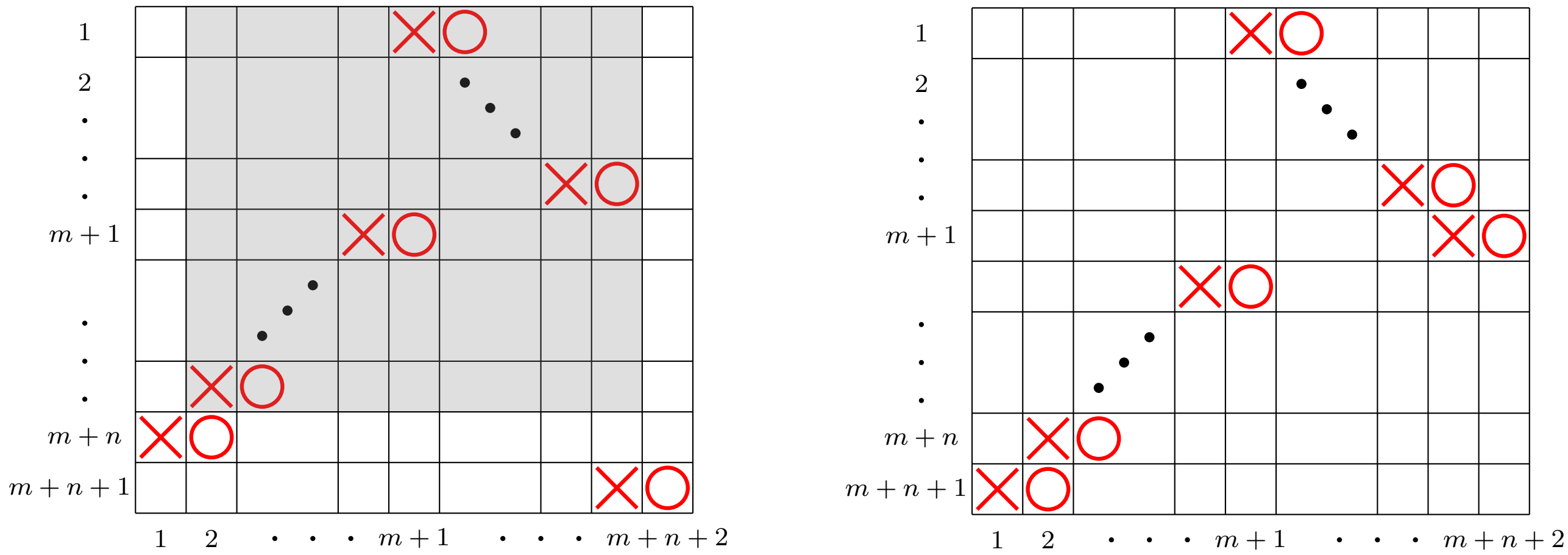

	\centering
	\begin{subfigure}{0.45\textwidth}
		\centering\tikzsetnextfilename{positive-inverted-staircase-step-3}
		\input{positive-inverted-staircase-step-3}
		\caption{Replacing $Z^\pm_{m-1,n-1}$ with $S^\pm_{n-1,m-1}$}
		\label{subfig:inverted-step-3}
	\end{subfigure}\hfill
	\begin{subfigure}{0.45\textwidth}
		\centering\tikzsetnextfilename{positive-inverted-staircase-step-4}
		\input{positive-inverted-staircase-step-4}
		\caption{Commutations produce $S^\pm_{n,m}$}
		\label{subfig:inverted-step-4}
	\end{subfigure}
	\caption{The inductive hypothesis in our proof of \Cref{prop:inverted-staircase} allows us to replace $Z^\pm_{m-1,n-1}$ with $S^\pm_{n-1,m-1}$, as on the left. Following this replacement, commuting the arc in row $m+n+1$ up to row $m+1$ completes the transformation of $Z^\pm_{m,n}$ into $S^\pm_{n,m}$.}
	\label{fig:staircase-spc2}
	\end{figure}
	
	After this sequence of SPC moves, rows $m+n$ and $m+n+1$ of $G$ each contain a length-$1$ arc, with that in row $m+n$ beginning in column $1$ and that in row $m+n+1$ beginning in column $m+n+1$. At the same time, the sub-grid of $G$ found in its first $m+n-1$ rows and columns $2,3,\ldots,m+n+1$ is identified with $Z^\pm_{m-1,n-1}$. According to our inductive hypothesis, $G$ is SPC-equivalent to the grid obtained by replacing this sub-grid with $S^\pm_{n-1,m-1}$. In fact, our inductive hypothesis ensures that this sub-grid replacement may be carried out via vertical commutations and bottom-to-top torus translations. Because $D$ contains no vertical arcs in the regions above and below this sub-grid of $G$, these SPC moves may be carried out on $D$, allowing us to replace the sub-grid $Z^\pm_{m-1,n-1}$ of $G$ with $S^\pm_{n-1,m-1}$ via vertical commutations and bottom-to-top torus translations on $D$. See \Cref{subfig:inverted-step-3}.
	
	Following this sub-grid replacement, we apply commutations to move the arc in row $m+n+1$ of $G$ up to row $m+1$ of $G$, as in \Cref{subfig:inverted-step-4}. After these commutations, we recognize $G$ as $S^\pm_{n,m}$, as desired. This concludes our inductive step, and thus establishes the result in the case where $i_j=0$ for $1\leq j\leq m+n+1$.

	Next, consider the case where $i_j=0$ for $1\leq j\leq m+n$, but $i_{m+n+1}$ may be nonzero. Under this assumption, the argument given above requires no modification --- because the arc of $D$ which is initially in row $m+1$ never leaves the sub-grid $G$, each of the vertical commutations applied during our argument remains valid when $i_{m+n+1}>0$.
	
	Finally, for a general multi-index $I$, observe that we may apply horizontal commutations to $Z^\pm_{m,n,I}$ to produce $Z^\pm_{m,n,I'}$, where
	\[
	I' = (0,0,\ldots,0,i_1+i_2+\cdots+i_{m+n+1}).
	\]
	That is, we apply horizontal commutations to move all empty columns to the right. If $D$ contains $Z^\pm_{m,n,I}$ as a sub-grid, then each of these commutations is also valid as a Cromwell move on $D$, so we may assume that $D$ instead contains $Z^\pm_{m,n,I'}$. After applying the SPC moves identified in the previous paragraph to replace $Z^\pm_{m,n,I'}$ with $S^\pm_{n,m,I'}$, we may return the empty columns to their original positions to produce $S^\pm_{n,m,I}$. That is, reversing the horizontal commutations which transformed $Z^\pm_{m,n,I}$ into $Z^\pm_{m,n,I'}$ will turn $S^\pm_{n,m,I'}$ into $S^\pm_{n,m,I}$. This concludes our proof.
	\end{proof}
	
	We conclude this section by observing that every standard-form purely dent-stabilized grid consists of a staircase and an inverted staircase.
	
	\begin{proposition}\label{prop:std-grids-from-staircases}
	For any integers $a,b,c,d\geq 0$, the sub-grids of $D_\pm(a,b,c,d)$ made up of its top $a+b+1$ rows and bottom $c+d+1$ rows, respectively, may be identified with $S^\pm_{a,b,I}$ and $Z^\mp_{c,d,J}$, respectively, where
	\[
	I = (i_1,i_2,\ldots,i_{a+b+1})
	\quad\text{and}\quad
	J = (j_1,j_2,\ldots,j_{c+d+1})
	\]
	satisfy $i_k=c\,\delta_{1k}+d\,\delta_{(a+b+1)k}$ and $j_{k}=(a+b)\,\delta_{(c+1)k}$, with $\delta_{ij}$ the Kronecker delta function.
	\end{proposition}

	The proof of \Cref{prop:std-grids-from-staircases} is a simple verification which we leave to the reader.
\section{SPC equivalence classes}\label{sec:spc-equivalence-classes}
	
	The goal of this section is to count the number of SPC- and unoriented SPC-equivalence classes. We begin by describing the relationship between the two orientations of a purely dent-stabilized grid in standard form. In particular, the following lemma shows that, for a standard grid diagram, if we both reverse the order of the four dent counts and reverse the orientation, then the SPC equivalence class is preserved.
	
	\begin{lemma}\label{lem:orientationReversalDent}
		Let $D(a,b,c,d)$ be a purely dent-stabilized grid in standard
		form. Then
		\[
		D_\pm(a,b,c,d)
		\sim_{\mathrm{SPC}}
		D_\mp(d,c,b,a).
		\]
	\end{lemma}
	
	\begin{figure}
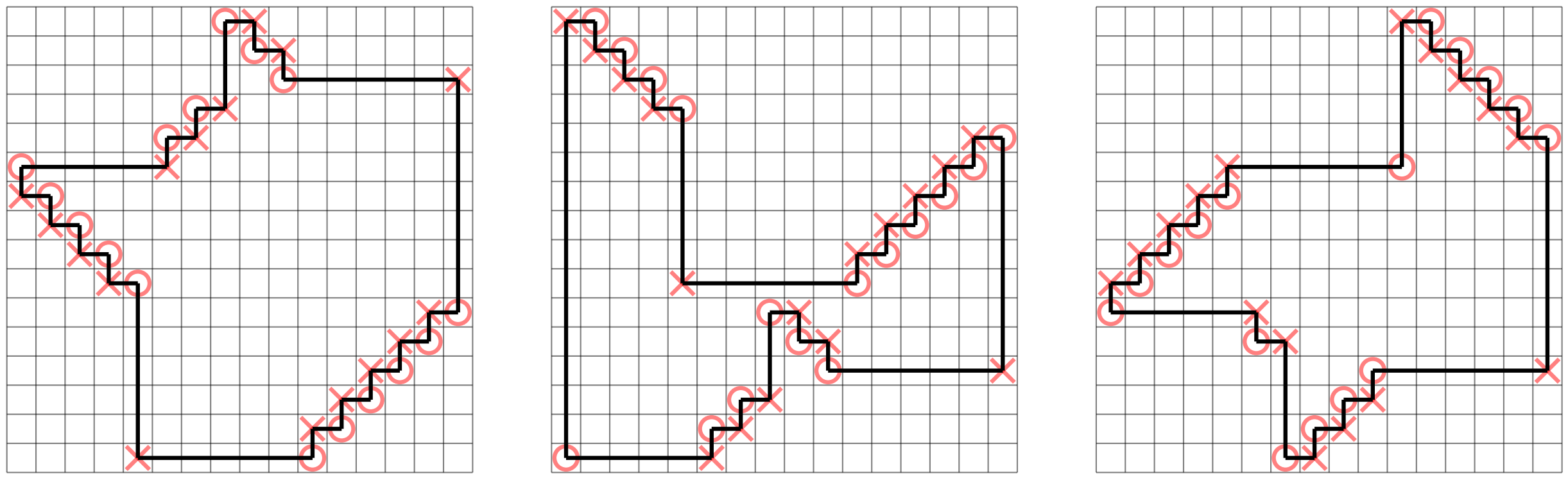

	\centering
		\begin{subfigure}{0.32\textwidth}
			\centering
			\input{reversal-lemma-start}
		\end{subfigure}\hfill
		\begin{subfigure}{0.32\textwidth}
			\centering
			\input{reversal-lemma-H}
		\end{subfigure}\hfill
		\begin{subfigure}{0.32\textwidth}
			\centering
			\input{reversal-lemma-end}
		\end{subfigure}
		\caption{The proof of \Cref{lem:orientationReversalDent} applies $c+d+1$ bottom-to-top torus translations to the grid diagram $D_\pm(a,b,c,d)$ to produce $H$, then applies \Cref{prop:inverted-staircase} to the top $c+d+1$ rows and, separately, to the bottom $a+b+1$ rows, to produce a crossingless grid diagram which is commutation-equivalent to $D_\mp(d,c,b,a)$. Here, this process is depicted for $D_-(3,2,4,5)$.}
		\label{fig:reversal-lemma}
	\end{figure}
	
	\begin{proof}
	Let $n = a+b+c+d+2$ be the grid number of $D_\pm(a,b,c,d)$. By \Cref{prop:std-grids-from-staircases}, we may identify the last $c+d+1$ rows of $D_\pm(a,b,c,d)$ with $Z^\mp_{c,d,I}$, where
	\[
	I = (i_1,i_2,\ldots,i_{c+d+1})
	\]
	is defined by $i_k=(a+b)\,\delta_{(c+1)k}$, for all $1\leq k\leq c+d+1$. Similarly, we may identify the first $a+b+1$ rows with $S^\pm_{a,b,J}$, where
	\[
	J = (j_1,j_2,\ldots,j_{a+b+1})
	\]
	is defined by $j_k=c\,\delta_{1k}+d\,\delta_{(a+b+1)k}$, for all $1\leq k\leq a+b+1$.
		
	Now, label the horizontal arcs of $D_\pm(a,b,c,d)$ from top to bottom by $R_1,R_2,\ldots,R_n$. Perform $c+d+1$ bottom-to-top torus translations and let the resulting grid be denoted by $H$. The rows of $H$, from top to bottom, therefore contain the arcs $R_{a+b+2},\ldots,R_{n-1},R_n,R_1,R_2,\ldots,R_{a+b+1}$, while the order of the vertical arcs is unchanged. The first $c+d+1$ rows of $H$ now form the sub-grid $Z^\mp_{c,d,I}$. By \Cref{prop:inverted-staircase}, we may replace this sub-grid with $S^\mp_{d,c,I}$. Denote the resulting grid by $H_0$.

	We next consider the final $a+b+1$ rows of $H_0$, which contain the arcs $R_1,R_2,\ldots,R_{a+b+1}$. The preceding staircase replacement changes only the specified $Z^\mp_{c,d,I}$ sub-grid, so these rows agree with those of $H$. The final $a+b+1$ rows thus form $S^\pm_{a,b,J}$. Applying \Cref{prop:inverted-staircase}, we replace this sub-grid $S^\pm_{a,b,J}$ with $Z^\pm_{b,a,J}$. Let $H_1$ denote the resulting grid.
		
	It remains to put the vertical arcs in their standard-form order. Consider the multi-indices
	\[
	I' = (i'_1,i'_2,\ldots,i'_{c+d+1})
	\quad\text{and}\quad
	J' = (j'_1,j'_2,\ldots,j'_{a+b+1})
	\]
	defined by
	\[
	i'_k := b\,\delta_{1k} + a\,\delta_{(c+d+1)k}
	\quad\text{and}\quad
	j'_k := (c+d)\,\delta_{(b+1)k}.
	\]
	Now the final $a+b+1$ rows of $H_1$ contain the sub-grid $Z^\pm_{b,a,J}$, and thus have $c$ empty columns immediately to the right of the leftmost column and $d$ empty columns immediately to the left of the rightmost column. Columns $c+2,c+3,\ldots,c+b+1$ of $H_1$ therefore contain vertical arcs which may be commuted to the left until they occupy columns $2,3,\ldots,b+1$. Likewise, the vertical arcs in columns $c+b+2,c+b+3,\ldots,c+b+a+1$ may be commuted to the right until they occupy columns $n-a,n-a+1,\ldots,n-1$. Following these commutations, one may verify that the bottommost $a+b+1$ rows of the grid diagram are identified with $Z^\pm_{b,a,J'}$, while the topmost $c+d+1$ rows are identified with $S^\mp_{d,c,I'}$. According to \Cref{prop:std-grids-from-staircases}, we have produced $D_\mp(d,c,b,a)$. Altogether, we see that
	\[
	D_\pm(a,b,c,d)
	\sim_{\mathrm{SPC}}
	D_\mp(d,c,b,a),
	\]
	as desired.		
	\end{proof}

	We can now classify exactly when two purely dent-stabilized grids in standard form are SPC- and unoriented SPC-equivalent.

	\begin{lemma}\label{lem:dentType}
	For any integers $a_j,b_j,c_j,d_j\geq 0$, $j=0,1$, the grid diagrams $D_\pm(a_0,b_0,c_0,d_0)$ and $D_\pm(a_1,b_1,c_1,d_1)$, with the same sign used in both, represent the same class:
	\begin{enumerate}
		\item in $G_n$ if and only if $(a_0,b_0,c_0,d_0)=(a_1,b_1,c_1,d_1)$.\label{part:oriented}
		\item in $\overline{G}_n$ if and only if $(a_0,b_0,c_0,d_0)=(a_1,b_1,c_1,d_1)$ or $(a_0,b_0,c_0,d_0)=(d_1,c_1,b_1,a_1)$.\label{part:unoriented}
	\end{enumerate}
	\end{lemma}
	
	\begin{proof}
	For any grid diagram $D$, let $D^\perp$ denote the grid diagram obtained by rotating $D$ clockwise through an angle of $90^\circ$. That is, $D^\perp$ contains an $X$ or $O$ in row $i$, column $j$ whenever $D$ contains an $X$ or $O$ in row $n+1-j$, column $i$.
	
	As discussed in \Cref{sec:background}, a grid diagram determines a Legendrian isotopy type, allowing us to speak sensibly of the \emph{Thurston-Bennequin} and \emph{rotation} numbers of $D$ and $D^\perp$. Because SPC moves preserve the Legendrian isotopy type represented by a grid diagram, the quadruple
	\[
	\mathcal{I}(D)=\left(\tb(D),\rot(D),\tb(D^\perp),\rot(D^\perp)\right)
	\]
	is invariant under SPC moves. Direct computation reveals that
	\[
	\mathcal{I}(D_\pm(a,b,c,d)) = (-1-a-d,\pm a \mp d, -1-b-c, \mp c\pm b),
	\]
	for any integers $a,b,c,d\geq 0$. In particular, a prescribed value $(t,r,t^\perp,r^\perp)$ for $\mathcal{I}(D_\pm(a,b,c,d))$ requires that
	\begin{equation}\label{eq:dents-from-classical-invariants}
	\begin{aligned}
		a&=\frac{-t-1\pm r}{2},
		&
		b&=\frac{-t^{\perp}-1\pm r^{\perp}}{2},\\
		c&=\frac{-t^{\perp}-1\mp r^{\perp}}{2},
		&
		d&=\frac{-t-1\mp r}{2}.
	\end{aligned}
	\end{equation}
	Now an SPC-equivalence between $D_\pm(a_0,b_0,c_0,d_0)$ and $D_\pm(a_1,b_1,c_1,d_1)$ would imply that
	\[
	\mathcal{I}(D_\pm(a_0,b_0,c_0,d_0)) = \mathcal{I}(D_\pm(a_1,b_1,c_1,d_1)),
	\]
	and from \Cref{eq:dents-from-classical-invariants} we would conclude that $(a_0,b_0,c_0,d_0)=(a_1,b_1,c_1,d_1)$. That is, the grid diagrams $D_\pm(a_0,b_0,c_0,d_0)$ and $D_\pm(a_1,b_1,c_1,d_1)$ are SPC-equivalent if and only if $(a_0,b_0,c_0,d_0)=(a_1,b_1,c_1,d_1)$, establishing \Cref{part:oriented}.
	
	Note that \Cref{lem:orientationReversalDent} immediately provides one direction of \Cref{part:unoriented}. Namely, if $(a_0,b_0,c_0,d_0)=(d_1,c_1,b_1,a_1)$, then \Cref{lem:orientationReversalDent} ensures that $D_\pm(a_0,b_0,c_0,d_0)$ and $D_\mp(a_1,b_1,c_1,d_1)$ represent the same equivalence class in $G_n$, and thus that $D_\pm(a_0,b_0,c_0,d_0)$ and $D_\pm(a_1,b_1,c_1,d_1)$ are equivalent in $\overline{G}_n$. Conversely, if $D_\pm(a_0,b_0,c_0,d_0)$ and $D_\pm(a_1,b_1,c_1,d_1)$ represent the same class in $\overline{G}_n$, then either
	\[
	D_\pm(a_0,b_0,c_0,d_0)\sim_{\mathrm{SPC}} D_\pm(a_1,b_1,c_1,d_1)
	\quad\text{or}\quad
	D_\pm(a_0,b_0,c_0,d_0)\sim_{\mathrm{SPC}} D_\mp(a_1,b_1,c_1,d_1).
	\]
	In the former case, \Cref{part:oriented} tells us that $(a_0,b_0,c_0,d_0)=(a_1,b_1,c_1,d_1)$; in the latter case, we invoke \Cref{lem:orientationReversalDent} to see that
	\[
	D_\pm(a_0,b_0,c_0,d_0) \sim_{\mathrm{SPC}} D_\mp(a_1,b_1,c_1,d_1) \sim_{\mathrm{SPC}} D_\pm(d_1,c_1,b_1,a_1)
	\]
	and then apply \Cref{part:oriented} to reach the conclusion $(a_0,b_0,c_0,d_0)=(d_1,c_1,b_1,a_1)$. So $D_\pm(a_0,b_0,c_0,d_0)$ and $D_\pm(a_1,b_1,c_1,d_1)$ represent the same equivalence class in $\overline{G}_n$ if and only if $(a_0,b_0,c_0,d_0)=(a_1,b_1,c_1,d_1)$ or $(a_0,b_0,c_0,d_0)=(d_1,c_1,b_1,a_1)$, as required by \Cref{part:unoriented}.
	\end{proof}

	We are now prepared to compute the cardinalities $|G_n|$ and $|\overline{G}_n|$.
	
	\begin{proof}[Proof of \Cref{thm:spc-classes}]
	According to \Cref{thm:strongtheorem1}, each equivalence class in $G_n$ and $\overline{G}_n$ is represented by a purely dent-stabilized grid diagram. Because every purely dent-stabilized grid diagram may be put into standard form via SPC moves, we conclude that each equivalence class in $G_n$ and $\overline{G}_n$ has a representative of the form $D_\pm(a,b,c,d)$, for some integers $a,b,c,d\geq 0$ with $a+b+c+d=n-2$. In fact, by using \Cref{lem:orientationReversalDent} to reverse the orientation, if necessary, we see that each equivalence class in $G_n$ and $\overline{G}_n$ has a representative of the form $D_+(a,b,c,d)$, for some integers $a,b,c,d\geq 0$ with $a+b+c+d=n-2$. \Cref{lem:dentType} then tells us that $G_n$ and $\overline{G}_n$ are in bijective correspondence with
	\[
	T_n := \{(a,b,c,d)\in\mathbb{Z}^4_{\geq 0}~|~a+b+c+d=n-2\}
	\quad\text{and}\quad
	\overline{T}_n := T_n/(a,b,c,d)\sim(d,c,b,a),
	\]
	respectively.
	
	The elements of $T_n$ are precisely the weak compositions of the integer $n-2$ into $4$ parts, and thus we have
	\[
	|G_n| = |T_n| = \dbinom{n+1}{3} = \dfrac{n(n^2-1)}{6},
	\]
	as desired. To determine the size of $\overline{T}_n$, consider the reversal action of $\mathbb{Z}_2=\{1,\rho\}$ on $T_n$ --- that is, the action defined by $\rho(a,b,c,d)=(d,c,b,a)$. Now $(a,b,c,d)$ is fixed by $\rho$ if and only if $a=d$ and $b=c$, so
	\[
	\Fix(\rho) = \{(a,b,b,a)\,|\,2a+2b=n-2\}.
	\]
	If $n$ is odd, then $\Fix(\rho)$ is evidently empty. If $n$ is even, then the elements of $\Fix(\rho)$ correspond to weak compositions of the integer $\tfrac{1}{2}n-1$ into $2$ parts, meaning that $|\Fix(\rho)|=\tfrac{1}{2}n$. Using Burnside's lemma, we then find that
	\[
	|\overline{G}_n| = \dfrac{1}{2}\left(|G_n|+0\right) = \dfrac{n(n^2-1)}{12}
	\]
	if $n$ is odd and
	\[
	|\overline{G}_n| = \dfrac{1}{2}\left(|G_n|+\frac{n}{2}\right) = \dfrac{n(n^2+2)}{12}
	\]
	if $n$ is even. This establishes \Cref{thm:spc-classes}.
	\end{proof}
\section{Cromwell hardness}\label{sec:applications}
	As discussed in the introduction, our main results are useful for the purposes of generating all grid diagrams for the unknot with grid number $n$. This collection of grid diagrams, once obtained, may then be searched for grid diagrams with special properties. In this final section, we search for grid diagrams which fail to be simplified to the grid diagram of grid number $2$ in a complexity-nonincreasing manner.

	\subsection{Generating grid diagrams for the unknot}
	When discussing the applications of our results in the introduction, we asserted that \Cref{thm:spc-classes} can help enumerate all elements of $\mathcal{G}_n$ with $n \geq 2$. In fact, we will enumerate the elements of each SPC-equivalence class of $\mathcal{G}_n$ separately; combining these $|G_n|=\tfrac{1}{6}n(n^2-1)$ lists produces an enumeration of the elements of $\mathcal{G}_n$. The following \emph{unknot generation algorithm} implements this strategy for a fixed grid size $n \geq 2$.
	
	We begin with an enumeration of the set
	\[
	\mathcal{S} := \{D_+(a,b,c,d)\,|\, a+b+c+d=n-2\}
	\]
	of \emph{seeds} as $S_1,S_2,\ldots,S_{|G_n|}$. The following algorithm produces all grid diagrams which are SPC-equivalent to the seed $S_k$, for a fixed $1\leq k\leq |G_n|$:
	\begin{enumerate}[label=(\arabic*)]
		\item Initialize the list $L_k$ as $L_k:=(S_k)$ and the integer $j$ as $j:=1$.
		
		\item If $j$ exceeds the length of $L_k$, the algorithm terminates. Otherwise, define $G$ to be the $j^{\text{th}}$ grid diagram in the list $L_k$ and define $P$ by
		\[
		P := \{\tilde{G}\in\mathcal{G}_n \,|\, \tilde{G} \text{~is obtained from~} G \text{~by a single SPC move}\}.
		\]
		Now append to $L_k$ those elements of $P$ which have not yet appeared in $L_k$, increase the value of $j$ by $1$, and repeat this step.
	\end{enumerate}
	Because there are finitely many grid diagrams of grid number $n$, this algorithm is guaranteed to terminate. Applying the algorithm to each of the seeds $S_1,S_2,\ldots,S_{|G_n|}$ produces a collection of lists $L_1,L_2,\ldots,L_{|G_n|}$, and, because every element of $\mathcal{G}_n$ is SPC-equivalent to some $D_+(a,b,c,d)$, these lists enumerate the elements of $\mathcal{G}_n$.

	\subsection{Searching for monotonic simplifications}
	Recall that we say that a grid diagram $D$ for the unknot is \emph{weakly Cromwell hard} if there is no sequence of Cromwell moves transforming $D$ into a grid diagram with grid number $2$ which is nonincreasing on both grid number and crossing number. This motivates the following definition:
	
	\begin{definition}
	The \emph{complexity} of a grid diagram $D$ for the unknot is the pair $\mathcal{C}(D) := (n,c)$, where $n$ is the grid number of $D$ and $c$ is the number of crossings of $D$. We define a partial order $\preceq$ on $\mathbb{N}\times\mathbb{N}$ by
	\[
	(n_1,c_1) \preceq (n_2,c_2)
	\iff
	n_1 \leq n_2 \text{~and~} c_1 \leq c_2.
	\]
	\end{definition}

	With this notion of complexity, notice that a grid diagram $D$ for the unknot is weakly Cromwell hard if for every sequence $D = D_1 \to \dots \to D_m = D_\pm(0,0,0,0),$ there exists an $i \in \{1, \ldots, m-1\}$ such that $\mathcal{C}(D_{i+1}) \npreceq \mathcal{C}(D_{i})$. That is, $D$ is weakly Cromwell hard if $D$ cannot be \emph{monotonically simplified} to $D_\pm(0,0,0,0)$ with respect to our notion of complexity.
		
	The following \emph{monotonic simplification check} algorithm allows us to determine whether a given grid diagram $D$ representing the unknot is weakly Cromwell hard:
	
	\begin{enumerate}[label=(\arabic*)]
		\item Initialize the list $L$ as $L:=(D)$ and the integer $j$ as $j:=1$.
		
		\item If $j$ exceeds the length of $L$, the algorithm terminates. Otherwise, define $G$ to be the $j^{\text{th}}$ grid diagram in the list $L$ and define $P$ by
		\[
		P := \{\tilde{G} \,|\, \tilde{G} \text{~is obtained from~} G \text{~by a single Cromwell move and~} \mathcal{C}(\tilde{G})\preceq \mathcal{C}(G)\}.
		\]
		Now append to $L$ those elements of $P$ which have not yet appeared in $L$. If $D_\pm(0,0,0,0) \in L$ after appending $P$, the algorithm terminates. Otherwise, increase the value of $j$ by $1$ and repeat this step.
	\end{enumerate}
	
	The algorithm must terminate, as no grid diagram of grid number greater than that of $D$ is considered.
	
	From the above formulation of weak Cromwell hardness, several symmetries among weakly Cromwell hard grid diagrams emerge. For instance, the complexity $\mathcal{C}$ is invariant under orientation reversal, and thus the question of weak Cromwell hardness depends only on the arc image of a grid diagram. We also see that complexity is invariant under rotation by $90^\circ$ and reflection across either the horizontal or vertical axis (i.e., reversing the orders of the rows or columns, respectively), as neither operation changes the grid or crossing number. So if $D$ is a weakly Cromwell hard  grid diagram, then we obtain up to $16$ weakly Cromwell hard grid diagrams: $8$ come from applying the elements of the dihedral group $D_4$ to $D$, and reversing the orientation of each of these produces another $8$. Therefore, when searching for Cromwell hard grid diagrams of the unknot, we do not need to search every SPC equivalence class. In particular, every SPC equivalence class has the form $[D_+(a,b,c,d)]$, for some integers $a,b,c,d\geq 0$ with $a+b+c+d=n-2$. Because $D_+(a,b,c,d)$ and $D_-(a,b,c,d)\sim_{\mathrm{SPC}}D_+(d,c,b,a)$ have the same arc image, it suffices to restrict attention to those SPC equivalence classes $[D_+(a,b,c,d)]$ which satisfy $a\leq d$. Reversing the rows of a grid diagram in the class $[D_+(a,b,c,d)]$ produces a grid diagram in $[D_+(b,a,d,c)]$, and thus we may further restrict attention to those SPC equivalence classes which satisfy $a\leq b$.
		
	Beyond these symmetries, other optimizations are available for the monotonic simplification check. The first is inductive on grid number. If we are checking a grid diagram $D$ of grid number $n$ and previous searching showed that there does not exist a weakly hard grid diagram with grid number less than $n$, then the appearance in the list $L$ of any grid diagram with grid number $n-1$ establishes that $D$ is monotonically simplifiable. The second optimization is similarly inductive on crossing number: if $\mathcal{C}(D)=(n,c)$ and all grid diagrams of grid number $n$ with fewer than $c$ crossings are known to be monotonically simplifiable, then the appearance in $L$ of any grid diagram with crossing number less than $c$ establishes the monotonic simplifiability of $D$. In either case, it is possible to conclude that a grid diagram is monotonically simplifiable before $D_\pm(0,0,0,0)$ has appeared in the list $L$.
	
	Once a weakly Cromwell hard grid diagram for the unknot has been found, modified versions of these optimizations remain in effect. For instance, if we find a weakly Cromwell hard grid diagram with complexity $(n_0,c_0)$, then a grid diagram $D$ is not verified to be monotonically simplifiable once a grid diagram of complexity $(n_0,c_0)$ appears in $L$. However, the appearance in $L$ of a grid diagram with complexity $(n,c)$ \emph{does} verify the monotonic simplifiability of $D$, provided every grid diagram of complexity $(n,c)$ is already known to be monotonically simplifiable.

	\subsection{The smallest weakly Cromwell hard grid diagrams}
	Using the algorithm described above, we have established the monotonic simplifiability of every grid diagram of the unknot of grid number at most $7$. The algorithm also produces a monotonic simplification for $165,872,992$ of the $165,873,120$ grid diagrams of the unknot having grid number $8$. For the remaining $128$ grid diagrams, the algorithm terminates without producing a monotonic simplification, certifying that these grid diagrams are weakly Cromwell hard and establishing \Cref{thm:min-weakly-hard}. The \emph{Mathematica} files used to produce these computation are available on the first author's website, linked at the end of this paper.
	
	As discussed above, there are certain symmetries among weakly Cromwell hard grid diagrams which lead to redundancies among this list of $128$ grid diagrams. The property of being weakly Cromwell hard is insensitive to orientation, and thus these $128$ grids represent just $64$ arc images. Moreover, rotating a weakly Cromwell hard grid diagram by $90^\circ$ evidently produces another weakly Cromwell hard grid diagram, and thus, using rotations, these $64$ arc images are easily produced from a collection of just $16$ arc images. Finally, reversing the vertical order of the rows of a weakly Cromwell hard grid diagram (equivalently, reflecting the grid diagram across a horizontal line through its center) produces another weakly Cromwell hard grid diagram, and thus $8$ arc images suffice to produce the full collection of $128$.
	
	\begin{figure}[ht]
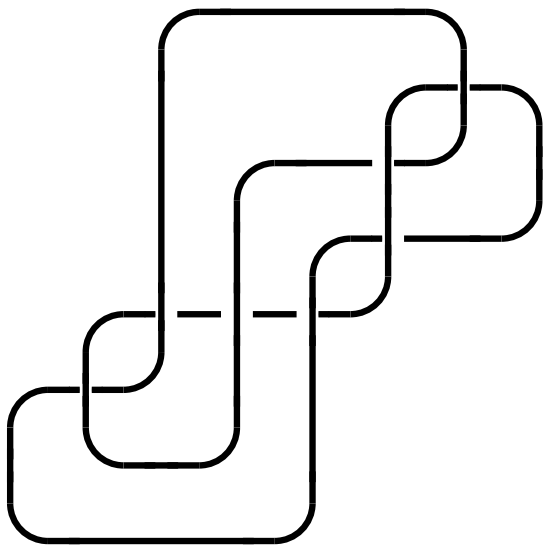
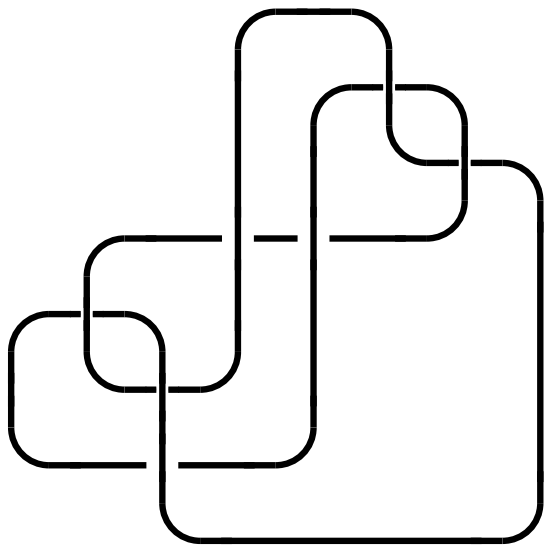
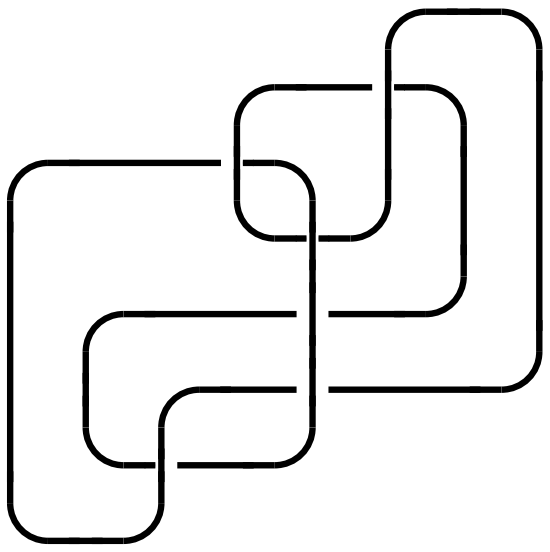
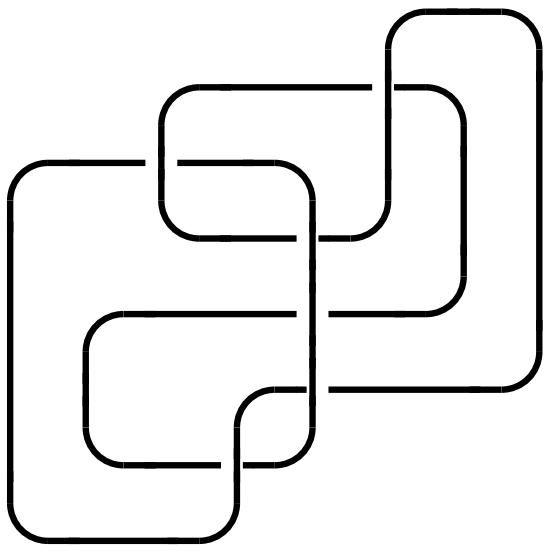
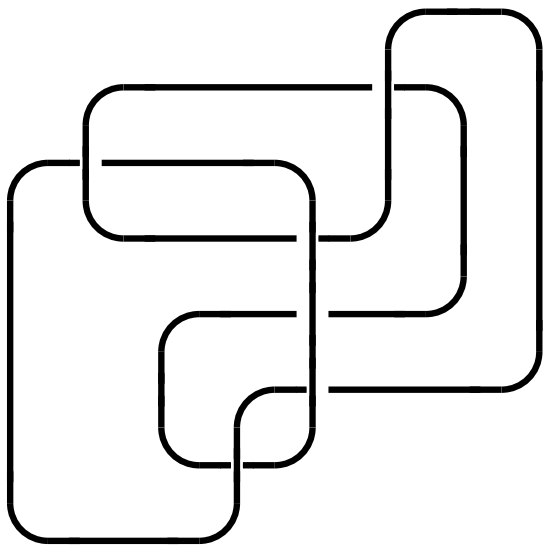

		\centering
		\begin{subfigure}{0.25\textwidth}
			\centering
			\input{hard-a.tex}
			\caption{Arc image $A$}
		\end{subfigure}\hfill
		\begin{subfigure}{0.25\textwidth}
			\centering
			\input{hard-b.tex}
			\caption{Arc image $B$}
		\end{subfigure}\hfill
		\begin{subfigure}{0.25\textwidth}
			\centering
			\input{hard-c.tex}
			\caption{Arc image $C$}
		\end{subfigure}\hfill
		\begin{subfigure}{0.25\textwidth}
			\centering
			\input{hard-d.tex}
			\caption{Arc image $D$}
		\end{subfigure}\hfill
		\begin{subfigure}{0.25\textwidth}
			\centering
			\input{hard-e.tex}
			\caption{Arc image $E$}
		\end{subfigure}\hfill
		\begin{subfigure}{0.25\textwidth}
			\centering
			\input{hard-f.tex}
			\caption{Arc image $F$}
		\end{subfigure}\hfill
		\begin{subfigure}{0.25\textwidth}
			\centering
			\input{hard-g.tex}
			\caption{Arc image $G$}
		\end{subfigure}\hfill
		\begin{subfigure}{0.25\textwidth}
			\centering
			\input{hard-h.tex}
			\caption{Arc image $H$}
		\end{subfigure}
		\caption{Representatives of weakly Cromwell hard grid diagrams}
		\label{fig:hard-grids}
	\end{figure}
	
	\Cref{fig:hard-grids} depicts $8$ arc images from which the $128$ grid diagrams may be generated via the rotations, reflections, and orientation reversals described above\footnote{Of course there is nothing canonical about these arc images; each is simply a representative of an equivalence class of hard grid diagrams.}. In fact, we notice that arc images $F$, $G$, and $H$ are related by vertical commutations (of their second, third, and fourth columns) which do not change the planar isotopy type of the diagram, making arc images $A$ through $F$ sufficient to generate all weakly hard grid diagrams of grid number $8$.
	
	\begin{table}
		\centering
		\begin{tblr}{colspec={cccc}}
	SPC class & $\mathbb{X}$-permutation & $\mathbb{O}$-permutation & arc image\\
	\hline
	$[D_+(0,2,3,1)]$ & $(3,4,8,7,6,1,2,5)$ & $(8,7,6,5,2,3,4,1)$ & \SetCell[r=2]{c}\input{hard-a-class-1}\\
	$[D_+(1,3,2,0)]$ & $(4,7,8,3,2,1,5,6)$ & $(8,5,6,7,4,3,2,1)$ & \\
	 \hline
	$[D_+(2,1,0,3)]$ & $(3,2,8,7,1,4,5,6)$ & $(1,4,3,2,5,6,7,8)$ & \SetCell[r=2]{c}\input{hard-a-class-2}\\
	$[D_+(3,0,1,2)]$ & $(3,4,5,8,2,1,7,6)$ & $(1,2,3,4,7,6,5,8)$ & \\
	 \hline
	$[D_+(2,0,1,3)]$ & $(6,5,1,2,3,8,7,4)$ & $(1,2,3,4,7,6,5,8)$ & \SetCell[r=2]{c}\input{hard-a-class-3}\\
	$[D_+(3,1,0,2)]$ & $(5,2,1,6,7,8,4,3)$ & $(1,4,3,2,5,6,7,8)$ & \\
	 \hline
	$[D_+(1,2,3,0)]$ & $(6,7,1,2,8,5,4,3)$ & $(8,5,6,7,4,3,2,1)$ & \SetCell[r=2]{c}\input{hard-a-class-4}\\
	$[D_+(0,3,2,1)]$ & $(6,5,4,1,7,8,2,3)$ & $(8,7,6,5,2,3,4,1)$ &  
\end{tblr}
		\caption{Weakly Cromwell hard grid diagrams produced by rotating and reflecting arc image $A$.}
		\label{table:hard-a-children}
	\end{table}
	
	The $128$ weakly Cromwell hard grid diagrams of grid number $8$ may be recovered explicitly from \Cref{table:hard-a-children,table:hard-b-children,table:hard-c-children,table:hard-d-children,table:hard-e-children,table:hard-f-children}. As done by other authors (e.g., \cite{ozsvath2015grid,petkova2026legendrian}), these tables present each grid diagram as a pair of permutations $\sigma_{\mathbb{X}},\sigma_{\mathbb{O}}\in S_n$. By definition, $\sigma_{\mathbb{X}}(i)=j$ tells us that the grid diagram contains an $X$ in row $n+1-j$, column $i$, and the permutation $\sigma_{\mathbb{O}}$ is defined analogously\footnote{While \cite{ozsvath2015grid} and \cite{petkova2026legendrian} use distinct conventions for $\sigma_{\mathbb{X}}$ and $\sigma_{\mathbb{O}}$, both number the rows of a grid diagram bottom-to-top. Our convention for the permutations matches that of \cite{ozsvath2015grid} once we recall our convention of numbering rows top-to-bottom.}. We present these permutations as $n$-tuples $(\sigma_{\mathbb{X}}(1),\ldots,\sigma_{\mathbb{X}}(n))$ and $(\sigma_{\mathbb{O}}(1),\ldots,\sigma_{\mathbb{O}}(n))$. While the tables list just $64$ grid diagrams, the orientation of each may be reversed by swapping the $\mathbb{X}$- and $\mathbb{O}$-permutations.
	
	\begin{table}
		\centering
		\begin{tblr}{colspec={cccc}}
	SPC class & $\mathbb{X}$-permutation & $\mathbb{O}$-permutation & arc image\\
	 \hline
	\SetCell[r=2]{c} $[D_+(1,3,1,1)]$ & $(1,5,4,3,2,6,7,8)$ & $(4,3,2,8,7,1,5,6)$ & \SetCell[r=4]{c} \input{hard-b-class-1}\\
	 & $(1,2,3,7,6,5,4,8)$ & $(3,4,8,2,1,7,6,5)$ & \\
	\SetCell[r=2]{c} $[D_+(1,1,3,1)]$ & $(1,5,4,3,2,6,7,8)$ & $(6,3,2,1,7,8,5,4)$ & \\
	 & $(1,2,3,7,6,5,4,8)$ & $(5,4,1,2,8,7,6,3)$ & \\
	 \hline
	\SetCell[r=2]{c} $[D_+(3,1,1,1)]$ & $(8,4,5,6,7,3,2,1)$ & $(3,6,7,8,2,1,4,5)$ & \SetCell[r=4]{c} \input{hard-b-class-2}\\
	 & $(8,7,6,2,3,4,5,1)$ & $(4,5,8,7,1,2,3,6)$ & \\
	\SetCell[r=2]{c} $[D_+(1,1,1,3)]$ & $(8,4,5,6,7,3,2,1)$ & $(5,6,7,1,2,8,4,3)$ & \\
	 & $(8,7,6,2,3,4,5,1)$ & $(6,5,1,7,8,2,3,4)$ & 
\end{tblr}
		\caption{Weakly Cromwell hard grid diagrams produced by rotating and reflecting arc image $B$.}
		\label{table:hard-b-children}
	\end{table}
	
	\begin{table}
		\centering
		\begin{tblr}{colspec={cccc}}
	SPC class & $\mathbb{X}$-permutation & $\mathbb{O}$-permutation & arc image\\
	\hline
	$[D_+(0,2,3,1)]$ & $(3,4,8,7,1,2,6,5)$ & $(8,7,6,2,3,5,4,1)$ & \SetCell[r=2]{c}\tikzsetnextfilename{hard-c-class-1}\input{hard-c-class-1}\\
	$[D_+(1,3,2,0)]$ & $(4,3,7,8,2,1,5,6)$ & $(8,5,4,6,7,3,2,1)$ & \\
	 \hline
	$[D_+(2,1,0,3)]$ & $(4,3,8,7,1,2,5,6)$ & $(1,5,4,2,3,6,7,8)$ & \SetCell[r=2]{c}\tikzsetnextfilename{hard-c-class-2}\input{hard-c-class-2}\\
	$[D_+(3,0,1,2)]$ & $(3,4,7,8,2,1,6,5)$ & $(1,2,3,6,7,5,4,8)$ & \\
	 \hline
	$[D_+(2,0,1,3)]$ & $(6,5,1,2,8,7,3,4)$ & $(1,2,3,7,6,4,5,8)$ & \SetCell[r=2]{c}\tikzsetnextfilename{hard-c-class-3}\input{hard-c-class-3}\\
	$[D_+(3,1,0,2)]$ & $(5,6,2,1,7,8,4,3)$ & $(1,4,5,3,2,6,7,8)$ & \\
	 \hline
	$[D_+(1,2,3,0)]$ & $(5,6,1,2,8,7,4,3)$ & $(8,4,5,7,6,3,2,1)$ & \SetCell[r=2]{c}\tikzsetnextfilename{hard-c-class-4}\input{hard-c-class-4}\\
	$[D_+(0,3,2,1)]$ & $(6,5,2,1,7,8,3,4)$ & $(8,7,6,3,2,4,5,1)$ &  
\end{tblr}
		\caption{Weakly Cromwell hard grid diagrams produced by rotating and reflecting arc image $C$.}
    \label{table:hard-c-children}
	\end{table}
	
	The SPC equivalence class of each weakly Cromwell hard grid diagram is identified in \Cref{table:hard-a-children,table:hard-b-children,table:hard-c-children,table:hard-d-children,table:hard-e-children,table:hard-f-children} by naming the unique standard-form representative $D_+$ to which the grid diagram is SPC-equivalent. Horizontal lines in the tables separate distinct unoriented SPC equivalence classes, and the arc image shown in the fourth column is that of the first grid diagram listed in each unoriented SPC equivalence class. There are evident structural similarities between \Cref{table:hard-a-children,table:hard-c-children}; between \Cref{table:hard-b-children,table:hard-e-children}; and between \Cref{table:hard-d-children,table:hard-f-children}, respectively --- including the repetition of certain permutations --- which we invite the reader to explore.
	
	\begin{table}
		\centering
		\begin{tblr}{colspec={cccc}}
	SPC class & $\mathbb{X}$-permutation & $\mathbb{O}$-permutation & arc image\\
	 \hline
	\SetCell[r=4]{c} $[D_+(1,2,2,1)]$ & $(1,2,3,6,5,4,8,7)$ & $(3,4,8,2,1,7,6,5)$ & \SetCell[r=4]{c} \tikzsetnextfilename{hard-d-class-1}\input{hard-d-class-1}\\
	 & $(2,1,5,4,3,6,7,8)$ & $(4,3,2,8,7,1,5,6)$ & \\
	 & $(1,2,3,6,5,4,8,7)$ & $(5,4,1,2,8,7,6,3)$ & \\
	 & $(2,1,5,4,3,6,7,8)$ & $(6,3,2,1,7,8,5,4)$ & \\
	 \hline
	\SetCell[r=4]{c} $[D_+(2,1,1,2)]$ & $(8,7,6,3,4,5,1,2)$ & $(4,5,8,7,1,2,3,6)$ & \SetCell[r=4]{c} \tikzsetnextfilename{hard-d-class-2}\input{hard-d-class-2}\\
	 & $(7,8,4,5,6,3,2,1)$ & $(3,6,7,8,2,1,4,5)$ & \\
	 & $(8,7,6,3,4,5,1,2)$ & $(6,5,1,7,8,2,3,4)$ & \\
	 & $(7,8,4,5,6,3,2,1)$ & $(5,6,7,1,2,8,4,3)$ & 
\end{tblr}
		\caption{Weakly Cromwell hard grid diagrams produced by rotating and reflecting arc image $D$.}
		\label{table:hard-d-children}
	\end{table}
	
	\begin{table}
		\centering
		\begin{tblr}{colspec={cccc}}
	SPC class & $\mathbb{X}$-permutation & $\mathbb{O}$-permutation & arc image\\
	 \hline
	\SetCell[r=2]{c} $[D_+(1,1,3,1)]$ & $(4,5,1,3,2,8,7,6)$ & $(2,3,4,8,7,6,5,1)$ & \SetCell[r=4]{c} \tikzsetnextfilename{hard-e-class-1}\input{hard-e-class-1}\\
	 & $(3,5,4,1,2,8,7,6)$ & $(8,1,2,3,7,6,5,4)$ & \\
	\SetCell[r=2]{c} $[D_+(1,3,1,1)]$ & $(3,2,1,7,6,8,4,5)$ & $(8,4,3,2,1,5,6,7)$ & \\
	 & $(3,2,1,7,8,5,4,6)$ & $(5,4,3,2,6,7,8,1)$ & \\
	 \hline
	\SetCell[r=2]{c} $[D_+(1,1,1,3)]$ & $(6,4,5,8,7,1,2,3)$ & $(1,8,7,6,2,3,4,5)$ & \SetCell[r=4]{c} \tikzsetnextfilename{hard-e-class-2}\input{hard-e-class-2}\\
	 & $(5,4,8,6,7,1,2,3)$ & $(7,6,5,1,2,3,4,8)$ & \\
	\SetCell[r=2]{c} $[D_+(3,1,1,1)]$ & $(6,7,8,2,1,4,5,3)$ & $(4,5,6,7,3,2,1,8)$ & \\
	 & $(6,7,8,2,3,1,5,4)$ & $(1,5,6,7,8,4,3,2)$ & 
\end{tblr}
		\caption{Weakly Cromwell hard grid diagrams produced by rotating and reflecting arc image $E$.}
		\label{table:hard-e-children}
	\end{table}
	
	We also notice that arc images $C$, $D$, and $E$ have $7$ crossings, while the other distinguished arc images all have $6$ crossings, leading us to the following question.
	\begin{question}
	Do there exist integers $\ell_n$ and $u_n$, for every integer $n\geq 8$, such that a weakly Cromwell hard grid diagram of grid number $n$ necessarily has crossing number between $\ell_n$ and $u_n$?
	\end{question}
	An affirmative answer to this question is not immediately unreasonable, as a low crossing number relative to the grid number suggests the availability of destabilizations, while a large crossing number suggests the availability of SPC moves which reduce the crossing number. In any case, the ability to identify precisely which crossing numbers are realizable via weakly Cromwell hard grid diagrams could be used to accelerate the search for weakly Cromwell hard grid diagrams of higher grid number.
		
	\begin{table}
		\centering
		\begin{tblr}{colspec={cccc}}
	SPC class & $\mathbb{X}$-permutation & $\mathbb{O}$-permutation & arc image\\
	 \hline
	\SetCell[r=12]{c} $[D_+(1,2,2,1)]$ & $(1,2,3,7,6,5,4,8)$ & $(6,4,1,5,2,8,7,3)$ & \SetCell[r=12]{c}\tikzsetnextfilename{hard-f-class-1}\input{hard-f-class-1}\\
	 & $(1,5,4,3,2,6,7,8)$ & $(6,2,1,7,4,8,5,3)$ & \\
	 & $(1,2,3,7,6,5,4,8)$ & $(3,5,8,2,4,1,7,6)$ & \\
	 & $(1,5,4,3,2,6,7,8)$ & $(3,2,8,5,7,1,4,6)$ & \\
	 & $(1,2,7,3,6,5,4,8)$ & $(6,4,5,1,2,8,7,3)$ & \\
	 & $(1,5,4,3,6,2,7,8)$ & $(6,2,1,7,8,4,5,3)$ & \\
	 & $(1,2,4,7,6,5,3,8)$ & $(4,5,8,2,3,1,7,6)$ & \\
	 & $(1,6,4,3,2,5,7,8)$ & $(3,2,8,6,7,1,4,5)$ & \\
	 & $(1,7,2,3,6,5,4,8)$ & $(6,5,4,1,2,8,7,3)$ & \\
	 & $(1,5,4,3,6,7,2,8)$ & $(6,2,1,7,8,5,4,3)$ & \\
	 & $(1,3,4,7,6,5,2,8)$ & $(4,5,8,3,2,1,7,6)$ & \\
	 & $(1,7,4,3,2,5,6,8)$ & $(3,2,8,7,6,1,4,5)$ & \\
	 \hline
	\SetCell[r=12]{c} $[D_+(2,1,1,2)]$ & $(8,7,6,2,3,4,5,1)$ & $(6,4,1,7,5,8,2,3)$ & \SetCell[r=12]{c}\tikzsetnextfilename{hard-f-class-2}\input{hard-f-class-2}\\
	 & $(8,4,5,6,7,3,2,1)$ & $(6,7,1,4,2,8,5,3)$ & \\
	 & $(8,7,6,2,3,4,5,1)$ & $(3,5,8,4,7,1,2,6)$ & \\
	 & $(8,4,5,6,7,3,2,1)$ & $(3,7,8,2,5,1,4,6)$ & \\
	 & $(8,7,5,2,3,4,6,1)$ & $(5,4,1,7,6,8,2,3)$ & \\
	 & $(8,3,5,6,7,4,2,1)$ & $(6,7,1,3,2,8,5,4)$ & \\
	 & $(8,7,2,6,3,4,5,1)$ & $(3,5,4,8,7,1,2,6)$ & \\
	 & $(8,4,5,6,3,7,2,1)$ & $(3,7,8,2,1,5,4,6)$ & \\
	 & $(8,6,5,2,3,4,7,1)$ & $(5,4,1,6,7,8,2,3)$ & \\
	 & $(8,2,5,6,7,4,3,1)$ & $(6,7,1,2,3,8,5,4)$ & \\
	 & $(8,2,7,6,3,4,5,1)$ & $(3,4,5,8,7,1,2,6)$ & \\
	 & $(8,4,5,6,3,2,7,1)$ & $(3,7,8,2,1,4,5,6)$ & 
\end{tblr}
		\caption{Weakly Cromwell hard grid diagrams produced by rotating and reflecting arc images $F$, $G$, and $H$.}
		\label{table:hard-f-children}
	\end{table}
	
\begin{figure}[ht]
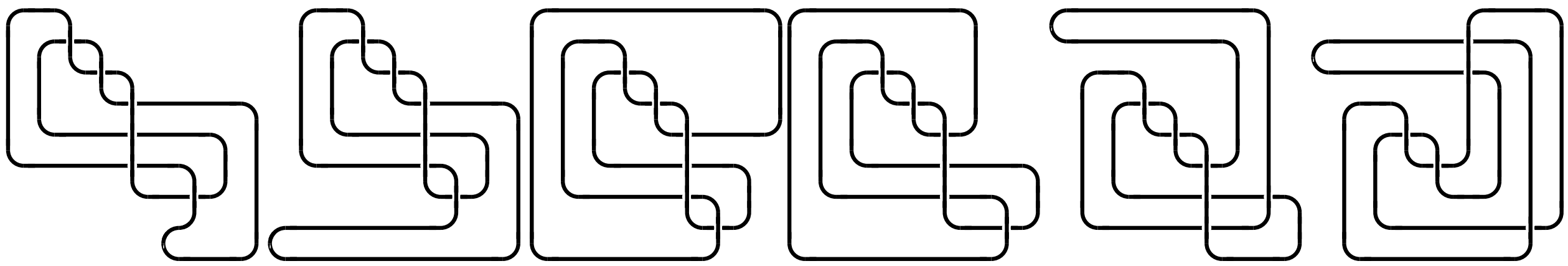

		\centering
		\begin{subfigure}{0.16\textwidth}
			\centering\tikzsetnextfilename{not-hard-movie-1}
			\input{not-hard-movie-1.tex}
		\end{subfigure}\hfill
		\begin{subfigure}{0.16\textwidth}
			\centering\tikzsetnextfilename{not-hard-movie-2}
			\input{not-hard-movie-2.tex}
		\end{subfigure}\hfill
		\begin{subfigure}{0.16\textwidth}
			\centering\tikzsetnextfilename{not-hard-movie-3}
			\input{not-hard-movie-3.tex}
		\end{subfigure}
		\begin{subfigure}{0.16\textwidth}
			\centering\tikzsetnextfilename{not-hard-movie-4}
			\input{not-hard-movie-4.tex}
		\end{subfigure}\hfill
		\begin{subfigure}{0.16\textwidth}
			\centering\tikzsetnextfilename{not-hard-movie-5}
			\input{not-hard-movie-5.tex}
		\end{subfigure}\hfill
		\begin{subfigure}{0.16\textwidth}
			\centering\tikzsetnextfilename{not-hard-movie-6}
			\input{not-hard-movie-6.tex}
		\end{subfigure}
		\caption{Stabilizing arc image $A$ once produces an arc image which is not weakly Cromwell hard. Here, the SPC moves are: a sequence of $5$ commutations; a torus translation; a pair of commutations; a torus translation; and another torus translation.}
		\label{fig:not-hard-movie}
	\end{figure}

	Finally, we note that none of the weakly Cromwell hard grid diagrams of grid number $8$ are strongly Cromwell hard --- that is, each admits a trivializing sequence of Cromwell moves which is nonincreasing on crossing number (though not on grid number). In particular, the reader may verify that each admits a trivializing sequence of Cromwell moves which includes just one stabilization. For instance, \Cref{fig:not-hard-movie} indicates how arc image $A$ may be simplified to a $5$-crossing grid diagram (which then admits a trivializing sequence which is nonincreasing on both crossing number and grid number) after a single stabilization. The following question thus remains:
	
	\begin{question}
	Do there exist any strongly Cromwell hard grid diagrams? If so, what is the smallest grid number of such a grid diagram?
	\end{question}
	
	\subsection{Weakly Cromwell hard grid diagrams of max-$\tb$ unknots}
	
	It was observed in \cite{breen2026hard} that many hard (in the Reidemeister sense) unknot diagrams cannot be realized as front diagrams for the unique Legendrian unknot in $(\mathbb{R}^3,\xi_{\mathrm{std}})$ with Thurston-Bennequin invariant $-1$. Nonetheless, \cite{breen2026hard} presents two types of hard Legendrian unknot: front diagrams for the max-$\tb$ unknot whose smooth resolutions are hard, as well as front diagrams for the max-$\tb$ unknot which are \emph{Legendrian hard}, meaning that the front diagram cannot be simplified to the $0$-crossing front diagram via a sequence of Legendrian Reidemeister moves which is nonincreasing on crossing number.
	
	Because grid diagrams determine Legendrian isotopy types, it is natural to ask in this setting questions analogous to those asked in \cite{breen2026hard}. We first observe, however, that the notions of weakly Cromwell hard and \emph{weakly Legendrian Cromwell hard} for a max-$\tb$ unknot are equivalent. Specifically, a sequence of Cromwell moves which is nonincreasing on the complexity of a grid diagram $D$ can only fail to preserve the Legendrian isotopy type of $D$ if it includes a destabilization which increases the Thurston-Bennequin invariant. If $\tb(D)=-1$, then no such destabilization can present itself, and thus a monotonic simplification of $D$ is immediately a Legendrian monotonic simplification.
	
	We will exhibit below a weakly Cromwell hard grid diagram for the unknot with $\tb(D)=-1$, certifying that such grid diagrams exist. Our work above tells us that any weakly Cromwell hard grid diagram for the max-$\tb$ unknot must have grid number at least $9$, since
	\[
	\tb(D_\pm(a,b,c,d)) = -1-a-d,
	\]
	and every weakly Cromwell hard grid diagram of grid number less than $9$ lies in an SPC equivalence class $[D_+(a,b,c,d)]$ with $a+d\geq 2$. This establishes \Cref{thm:min-weakly-hard-max-tb}.
	
\begin{figure}[t]
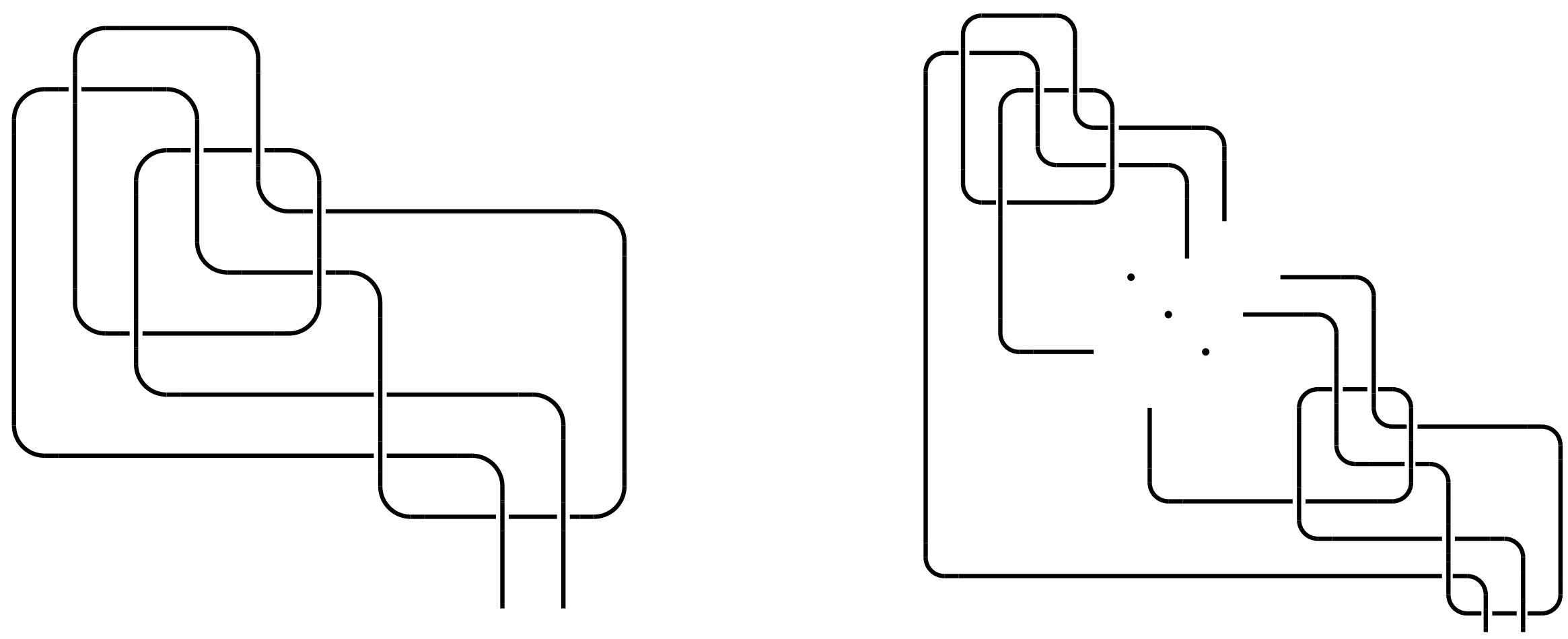

	    \centering
	    \begin{subfigure}{0.45\textwidth}
			\centering\tikzsetnextfilename{building-block}
		\input{building-block}
		\caption{A $10\times 10$ building block.}
    		\label{fig:building-block}
	    \end{subfigure}\hfill
	    \begin{subfigure}{0.45\textwidth}
			\centering\tikzsetnextfilename{big-building-block}
		\input{big-building-block}
		\caption{A $(6+4n)\times(6+4n)$ building block.}
    		\label{fig:big-building-block}
	    \end{subfigure}
	    \caption{Building blocks for weakly Cromwell hard grid diagrams of the max-$\tb$ unknot.}
	    \label{fig:building-blocks}
	\end{figure}

\begin{figure}[ht]
		\centering\tikzsetnextfilename{max-tb}
	    \input{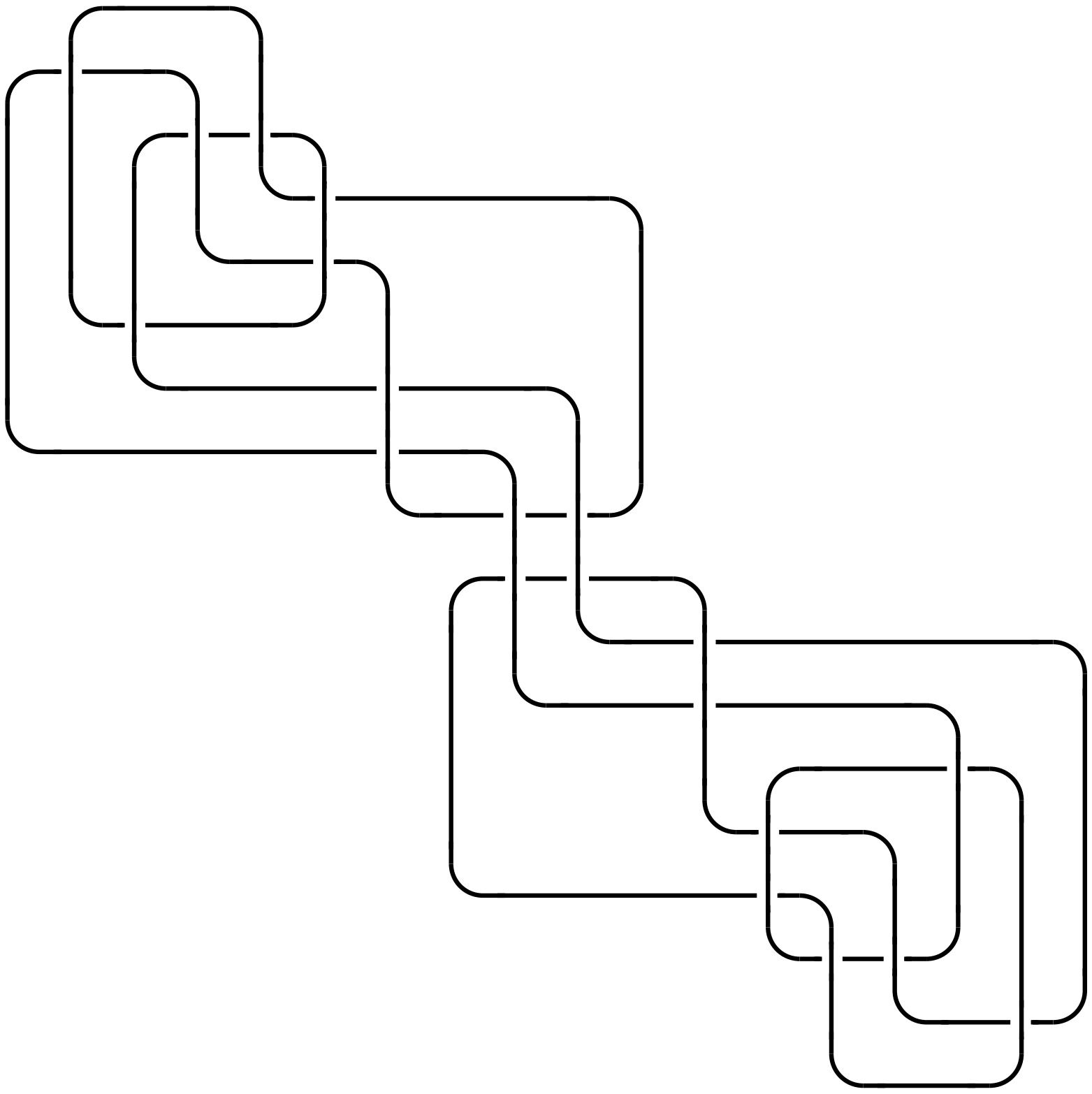}
	    \caption{The arc image of the grid diagram with $\mathbb{X}$-permutation $(2, 4, 3, 5, 7, 6, 10, 1, 11, 8, 18, 9, 13, 12, 14, 16, 15, 17)$ and $\mathbb{O}$-permutation $(5, 1, 6, 11, 4, 2, 3, 9, 7, 12, 10, 16, 17, 15, 8, 13, 18, 14)$. This is a weakly Cromwell hard grid diagram for the unknot with $\tb(D)=-1$ and grid number $18$.}
    	    \label{fig:max-tb}
	\end{figure}

	We now produce explicit weakly Cromwell hard grid diagrams for the max-$\tb$ unknot. Our construction is based on the building blocks introduced in \cite[Lemma 3.2]{breen2026hard}. In particular, \cite[Figure 29]{breen2026hard} presents two partial front diagrams, which may be labeled $B_0$ and $B_1$. A realization of $B_1$ as a partial arc image is seen in \Cref{fig:building-block}. By rotating $B_1$ through an angle of $180^\circ$ and joining the resulting partial arc image to $B_1$ in the obvious manner, we produce the arc image for the unknot seen in \Cref{fig:max-tb}. Application of the monotonic simplification algorithm allows us to conclude that this grid diagram is weakly Cromwell hard, and direct computation reveals that its Thurston-Bennequin invariant is $-1$. Note that in this case, the monotonic simplification algorithm can be carried out by hand: applying any torus translation to the arc image in \Cref{fig:max-tb} would increase the crossing number; there are two (vertical) arcs which admit crossing-nonincreasing commutations --- leading to two more available torus translations --- but these torus translations also increase the crossing number. We remark that applying the same construction to building block $B_0$ of \cite[Figure 29]{breen2026hard}, or to the pair $B_0$, $B_1$, produces a grid diagram which is not weakly Cromwell hard.
	
	As in \cite{breen2026hard}, the construction above extends to an infinite family. Mimicking \cite[Figure 32]{breen2026hard}, additional loops may be inserted into $B_1$ to produce the $(6+4n)\times(6+4n)$ partial arc image seen in \Cref{fig:big-building-block}, which we label $B_n$. By rotating $B_n$ through an angle of $180^\circ$ and joining with $B_1$ in the same manner as above, we obtain a weakly Cromwell hard grid diagram for the max-$\tb$ unknot of grid number $14+4n$. Weak Cromwell hardness is verified as in the $n=1$ case: each of the four torus translations would increase crossing number, and the two vertical arcs which admit crossing-nonincreasing commutations lead to torus translations which increase crossing number. In particular, it follows that there exist weakly Cromwell hard max-$\tb$ unknot grid diagrams of arbitrarily large grid number.

\bibliography{references}
\bibliographystyle{amsalpha}

\end{document}